\documentclass{amsart}
\usepackage{amsmath}
\usepackage{amssymb}
\usepackage{euscript}
\usepackage[all]{xy}
\usepackage{varioref}

\newtheorem{theorem}{Theorem}
\newtheorem{corollary}{Corollary}
\newtheorem{lemma}{Lemma}
\newtheorem{proposition}{Proposition}

\newtheorem{definition}{Definition}
\input cyracc.def

\newcommand{\Span}[1]{\ensuremath{ 
  \operatorname{span} 
   \left \{ {#1} \right \}
}}
\newcommand{\pd}[2]{\ensuremath{ \frac{\partial {#1}}{\partial {#2}} }}
\newcommand{\R}[1]{\ensuremath{\mathbb{R}^{#1}}}
\newcommand{\C}[1]{\ensuremath{\mathbb{C}^{#1}}}
\newcommand{\CP}[1]{\ensuremath{\mathbb{CP}^{#1}}}
\newcommand{\Gro}[2]{\widetilde{\operatorname{Gr}}_{ #1 } \left ( #2 \right ) }
\newcommand{\GL}[1]{ \operatorname{GL} \left ( {#1} \right ) }
\newcommand{\gl}[1]{ \operatorname{gl} \left ( {#1} \right ) }
\newcommand{\GLp}[1]{ \operatorname{GL}^+ \left ( {#1} \right ) }
\newcommand{\SL}[1]{ \operatorname{SL} \left ( {#1} \right ) }
\newcommand{\SO}[1]{ \operatorname{SO} \left ( {#1} \right ) }
\newcommand{\Lin}[2]{ \operatorname{Lin} \left ( {#1} , {#2} \right ) }

\newcommand{\Cstrucs}[1]{\mathcal{J} \left ( #1 \right )}
\newcommand{\Lm}[2]{\ensuremath{\Lambda^{#1} \left ( {#2} \right )}}
\newcommand{\nForms}[2]{\ensuremath{\Omega^{#1} \left ( {#2} \right )}}

\newcommand{\tr}{\ensuremath{\operatorname{tr}}}
\newcommand{\Ad}{\ensuremath{\operatorname{Ad}}}
\newcommand{\ModSp}[2]{\ensuremath{{#2}^{#1}}}
\begin{document}
\title{Dual Curves and Pseudoholomorphic Curves}
\author{Benjamin McKay}
\address{University of Utah}
\email{mckay@math.utah.edu}
\date{\today}
\begin{abstract} A notion of dual curve for
pseudoholomorphic curves in 4--manifolds 
turns out to be possible only if the notion of almost
complex structure structure is slightly generalized.
The resulting structure is as easy (perhaps
easier) to work with, and yields many
analogues of results in complex surface theory,
using a description of the local geometry 
via Cartan's method of equivalence.
Duality then uncovers a new
infinite-dimensional family of geometric
integrable systems. These are the first
steps toward geometry on moduli spaces 
of pseudoholomorphic curves.
\end{abstract}
\maketitle
\section{Introduction} 
Consider the notion of dual curves in the
projective plane: the space of lines
in the plane forms a plane called the
\emph{dual plane}, and the dual of a 
plane curve
is its collection of tangent lines, thought
of as a curve in the dual plane. Moreover,
the dual of the dual is the origin curve.
And the dual of a point is a line. Does
this picture work for almost complex structures
on the complex projective plane? Mikhael
Gromov (see \cite{Gromov:1985}) showed that
the dual plane, taken as the space of
``complex lines'', i.e.\ pseudoholomorphic spheres in the
homology class of a complex line, 
is a 4-manifold (if the almost complex
structure is tame), but is
it equipped with an almost complex structure?
I would like to use tangential approximation
by ``complex lines'' to define dual curves,
and would like them to be pseudoholomorphic.
Naive dimension count suggests that this
should determine the almost complex structure
on the dual projective plane. However,
we will see that there is no such almost
complex structure, except
in the case we have already considered:
the standard complex projective plane.
Except in that case, the dual projective plane is not 
almost complex in any manner which would
render  the dual curves pseudoholomorphic.

  There is a more fundamental question:
why think about almost complex structures in
the first place? We know that they 
exist on symplectic manifolds, and that they
have pseudoholomorphic curves which are
very similar to holomorphic
curves in complex manifolds.
The curves are usually used to probe symplectic
geometry. But what we use about them is
that they form a family of surfaces 
defined by local conditions and depending
on several functions of one real variable,
and have carefully controlled singularities.
In particular, they cannot crease. On top of
this, we can use a symplectic structure to ``tame''
these curves---to provide an a priori estimate
preventing all but finitely many singularities
of a mild nature.

  The largest family of determined first
order equations
on 4-manifolds with these same properties,
a family I will call
\emph{pseudocomplex structures},
are described in the same article of Gromov
\cite{Gromov:1985} (there called
\emph{elliptic systems}). I will
use them instead of almost complex
structures, and discover that my
construction of dual curves on a 
dual projective plane above will 
always work in
the category of pseudocomplex structures.
Moreover, the analysis of these equations, both
local and global,
is as easy as the analysis
of pseudoholomorphic curves for
almost complex structures.
It is the author's belief that the
pseudocomplex structure 
is more natural and
fundamental than the almost
complex structure.
To demonstrate this, proofs
are provided or sketched of 
analogues for the
basic results of the theory of complex 
surfaces.

  An application of this
idea will be given: an explicit geometric
description via dual curves of an infinite-dimensional
family of integrable systems of
pseudocomplex structures
first described indirectly by Gaston Darboux
in \cite{Darboux:1894}.

\subsection{Other work on this subject} 
Much of the material described here
began in the author's thesis
\cite{McKay:1999}. The interested
reader will find  more details
of calculations there which I have omitted
here for lack of space. The author wishes
to acknowledge that some of the results of this
paper are nearly identical to those
independently discovered by Jean-Claude
Sikorav in \cite{Sikorav:2000}, which appeared
in preprint form somewhat earlier
than this article, although later than
my doctoral thesis. Sikorav's approach
treats singularities more fully, 
elegantly and economically, by employing
results of Micallef \& White, while my
approach has its advantages in organizing
the calculations of local invariants,
and thereby proving Sikorav's Main
Conjecture. No new mathematical material has
been added to this article since the
appearance of Sikorav's preprint on
LANL.

\section{Linear algebra and microlocal geometry}

In this section we will first present
a heuristic discussion of ellipticity,
and then provide proofs.

\subsection{Rough ideas of ellipticity}
Consider a system of differential equations
for surfaces in a 4-manifold whose 
general solution depends on
arbitrary functions of one variable,
and which has no creasing of solutions,
i.e.\ in any metric, mean curvature
of solutions is bounded by a constant
on any compact set. Calculation shows
that if the system is real analytic,
then this system must be elliptic
and determined. If we insist further
that the system be first order,
we call such a system a
\emph{pseudocomplex structure} by analogy with the 
Cauchy--Riemann equations of complex
curves in a complex surface or
almost complex 4-manifold.
I will refer to the solutions of
a pseudocomplex structure
as its \emph{curves}.

  A more formal definition (exactly
as in \cite{Gromov:1985}): a
\emph{pseudocomplex structure} 
is a six-dimensional manifold
$E$ which is a fiber subbundle
of $\Gro{2}{TM} \to M$, the bundle
of oriented 2-planes in the tangent
spaces of a 4-manifold $M$,
such that the requirement that
a surface $\Sigma \subset M$
should have tangent spaces lying
in $E$ is equivalent to an elliptic
system of partial differential equations.
Call such a surface $\Sigma$ an
\emph{$E$ curve}.
If the fibers of $E \to M$ are compact,
then we will say that $E$ is 
\emph{proper} (because the map $E \to M$
is a proper map precisely when
the fibers are compact).
For example, take $M$ a complex
surface and let $E$ be the collection
of complex lines in the tangent spaces
of $M$. 

  In a single fiber, $E_m \subset \Gro{2}{T_m M}$
above a point $m \in M$, we have a surface
in the Grassmannian of 2-planes. Projectivizing,
this is a surface in the Grassmannian
of projective lines in projective 3-space,
in other words a \emph{line congruence}.
We will now study these in detail.

   Take $V$ any four-dimensional
vector space, and $X \to \Gro{2}{V}$
any immersed surface. We will call $X$
a line congruence.
A complex structure on $V$, say $J:V \to V$,
can be represented as a line congruence
by its Riemann sphere of complex lines,
which form a 2 parameter family of 
oriented real 2-planes. An arbitrary 
immersed surface $X \to \Gro{2}{V}$ 
will be called \emph{elliptic}
if it is tangent at each point to the
Riemann sphere of some complex structure.
Such a structure is never unique. This
notion of ellipticity is exactly the
one needed for our definition above
of pseudocomplex structure:
the fibers $E_m \to \Gro{2}{T_m M}$
of $E$ are elliptic line congruences
precisely when the differential equation 
for surfaces in $M$ which is determined by $E$
is elliptic.
Another definition of ellipticity:
thinking of a line congruence as a
2-parameter family $X$ of lines in 
$\mathbb{P}^3$, we want these
lines to twist about any transverse
line in $\mathbb{P}^3$ with nowhere
vanishing angular momentum. The
reader might want to draw a
picture of a Riemann sphere as
a line congruence to get a geometric
feeling for this. A compact elliptic
line congruence is precisely a fibration
of $\mathbb{P}^3$ by projective lines.

Yet another approach to define
ellipticity: in $\GL{V}$ invariant fashion
we can identify
\[
T_p \Gro{2}{V} \cong \Lin{p}{V/p}
\]
so that the determinant of these linear
maps, well defined up to scaling, provides
a conformal quadratic form on the tangent spaces
of the Grassmannian, of signature $(2,2):$
\[
\det
\begin{pmatrix}
p+q & r+s \\
r-s & p-q
\end{pmatrix}
= p^2 + s^2 - q^2 - r^2.
\] 
This conformal
quadratic form restricts to $X$ to be a positive
definite conformal structure precisely
when $X$ is an elliptic line congruence.
The stabilizer of a 2-plane $p \in \Gro{2}{V}$
acts transitively on positive definite
2-planes $\Pi \subset T_p \Gro{2}{V}$,
i.e.\ on positive definite 2-planes in
$\Lin{p}{V/p}$.  Therefore
any two elliptic line congruences can
be made to match up to first order
at any points by linear transformation.

\subsection{Making ellipticity more precise}

Let us now prove the statements we have
just made about ellipticity. We have 
defined ellipticity of a line congruence $X \to \Gro{2}{V}$
to be the property that its tangent spaces 
$T_x X \subset T_P \Gro{2}{V} = \Lin{P}{V/P}$
contain no linear maps of rank 1,
or equivalently that the conformal
quadratic form $\det$ be definite.
We can suppose by choice of orientation
of $V$ that it is positive definite.
For a pair of partial differential
equations
\[
F^i\left(x^1,x^2,u^1,u^2,\pd{u^j}{x^k}\right)=0, \ i=1,2,
\]
we wish to define ellipticity to be
the absence of real points in the characteristic
variety. It remains to show that these 
definitions coincide in the sense outlined
above.

\begin{lemma}
A pair of independent equations on the first
derivatives of two functions of two variables
is an elliptic system of partial differential
equations (i.e.\ there are no real points in
the characteristic variety of any integral element) 
precisely when the two equations
cut out, for any fixed values of independent
and dependent variables, a surface in the 
``space of first derivatives'' (i.e.\ the
Grassmannian of 2-planes) which is an elliptic
line congruence (i.e.\ the surface is 
definite for the quadratic form $\det$
on the Grassmannian).
\end{lemma}
\begin{proof}
First let us note that for this
pair of equations to be independent equations
on the first derivatives of the two functions
$u^1,u^2$ of the two variables $x^1,x^2$, we need
to ask that, if we write these two equations as
\[
F^i\left(x^1,x^2,u^1,u^2,p^j_k\right)=0
\]
with fictitious variables $p^j_k$ replacing
$\pd{u^j}{x^k}$, then
\[
\operatorname{rk}
\begin{pmatrix}
\pd{F^1}{p^1_1} & \pd{F^1}{p^1_2} & \pd{F^1}{p^2_1} &  \pd{F^1}{p^2_2} \\
\pd{F^2}{p^1_1} & \pd{F^2}{p^1_2} & \pd{F^2}{p^2_1} &  \pd{F^2}{p^2_2} \\
\end{pmatrix}
= 2.
\]
The characteristic variety is then defined by
the equations
\[
0 = \det
\begin{pmatrix}
\pd{F^1}{p^1_k} \xi^k & \pd{F^1}{p^2_k} \xi^k \\
\pd{F^2}{p^1_k} \xi^k & \pd{F^2}{p^2_k} \xi^k \\
\end{pmatrix}
\]
on some new variables $\xi_1,\xi_2$
(see Bryant et al.\ \cite{BCGGG:1991}).
We are asking that these equations have no
real solutions except $\xi=0.$ 
Write
\[
DF^i = 
\begin{pmatrix}
\pd{F^1}{p^1_1} & \pd{F^1}{p^2_1} \\
\pd{F^1}{p^1_2} & \pd{F^1}{p^2_1} 
\end{pmatrix}
\]
and write $\xi=\left(\xi_1, \xi_2\right)$
a row vector.
If we have such a solution $\xi,$ then
we obtain a linear relation between the
row vectors \( \xi DF^1 \) and \( \xi DF^2, \)
which we can suppose, switching subscripts
if needed, is \( \xi DF^1 = \lambda \xi DF^2. \)
But then \( DF^1 - \lambda DF^2 \) has rank one.
Conversely, if a linear combination 
of \( DF^1 \) and \( DF^2 \) has rank one,
then we can pick such a $\xi,$ and the
characteristic variety is not empty.

Thinking now of our differential equations
as defining a submanifold of the Grassmann
bundle, we think of the 2-planes 
\[
du^i = p^i_j \, dx^j
\]
as our family of 2-planes, with the equations
$F^1=F^2=0$ determining 2 constraints on
those 2-planes, i.e.\ as the
equations of our subbundle. 
Therefore the tangent space
to this family of 2-planes is given by
the equations $dF^1=dF^2=0.$ Taking the
Grassmannian of 2-planes at a fixed
point $\left(x^1,x^2,u^1,u^2\right)$,
the tangent plane to our family of 2-planes
is precisely cut out by the equations
$dF^1=dF^2=0$
with the added conditions that $dx^1=dx^2=du^1=du^2=0$
forcing us to stay at a particular point.
But then
\[
dF^1=\pd{F^1}{p^i_j} dp^i_j
\]
so that the equations on the tangent
space are precisely
\[
DF^1 \, dp = DF^2 \, dp = 0. 
\]

The identification of tangent vectors
$v \in T_P \Gro{2}{V}$ with linear
maps is carried out as follows: a curve
in $\Gro{2}{V}$ is a one parameter
family $P(t)$  of 2-planes: $P(t) \subset V.$
Taking any linear map $\phi(t) : V \to W$
with kernel $P(t)$, assuming that $\phi(t)$
depends smoothly on $t$, and letting 
$\left[\phi(t)\right] : V/P(t) \to W$
be the induced linear map on $V/P(t),$
we can identify 
the vector $P'(0) \in T_{P(0)} \Gro{2}{V}$
with the linear map
\[
\left[\phi(0)\right]^{-1} \phi'(0) : P(0) \to V/P(0).
\] 
We then have the conformal quadratic form
defined by taking determinant, which
depends on a choice of volume element
in $P(0)$ and in $V$, so that without
choosing a volume element, the quadratic
form is only defined up to conformal 
transformations.

Returning to our partial differential
equation, in our local coordinate formulation, we can
take $V$ to have coordinates $x^1,x^2,u^1,u^2$
and let
$P(0)=\left\{(x,0)|x \in \R{2}\right\}$
and take $\dot{p}=\left(\dot{p}^i_j\right)$
any tangent vector to our fiber, i.e.\ 
so that $DF^1 \dot{p} = DF^2 \dot{p} = 0,$
and find that using 
\[
\phi(t) : (x,u) \mapsto u+t \dot{p} x
\]
the matrix representing the
vector $\dot{p}$ is precisely $\dot{p}$.
Note that we have assumed that the
tangent vector is tangent to a point $(x,u,p)$
with $p=0;$ the requirement that our point have $p=0$ 
can be arranged by a change of coordinates,
or we can use $\phi(t) = (x,u) \mapsto u+\left(p+t\dot{p}\right)x.$

It is clear that a 2-plane
is positive definite precisely when its  dual
space is positive definite. Note that the dual
space of $\Lin{P}{V/P}$ is $\Lin{V/P}{P},$
and is equipped with the same quadratic form. 
We have seen that the absence of real 
points in the characteristic
variety is precisely the definiteness of the
tangent planes to the fibers of our subbundle. We can 
assume that they are positive definite
by choice of orientation of $V$.
\end{proof}

\subsection{Properties of line congruences}

The only known deep result on line congruences
was found by Gluck and Warner 
\cite{Gluck/Warner:1983}.
If we impose a positive definite 
conformal structure on $V$ and
an orientation, then we have a splitting
\[
\Lm{2}{V} = \Lm{2+}{V} \oplus \Lm{2-}{V}
\]
of 2-forms into self-dual and anti-self-dual 2-forms,
and the Pl\"ucker embedding 
\[
\Gro{2}{V} \hookrightarrow \Lm{2}{V} / \mathbb{R}^+
\]
which has as image the projectivization of the locus
\( \{ \xi | \xi^2 = 0 \}. \) This gives a
diffeomorphism
\[
\Gro{2}{V} \cong 
\Lm{2+}{V} / \mathbb{R}^+ \times \Lm{2-}{V} / \mathbb{R}^+
\cong S^{2+} \times S^{2-}
\]
into a product of spheres, which is invariant 
under the orientation
preserving conformal group of $V$.

\begin{theorem}[Gluck and Warner]
Given an immersed surface $\phi : X \subset \Gro{2}{V}$ (i.e.\
a line congruence) write
\[
\phi_+ : X \to S^{2+}
\]
and $\phi_-$ similarly, for the induced projections
to the two spheres $S^{2+}$ and $S^{2-}$.
A line congruence is elliptic
precisely when $|\phi'_-| > |\phi'_+|$
at all points of $X$.
\end{theorem} 
\begin{proof}
Indeed we will see that 
the $(2,2)$ conformal structure on $\Gro{2}{V}$
is just $| \cdot_{-} |^2 - | \cdot_{+} |^2$.
The result is surprising since the notion
of ellipticity of a line congruence is
invariant under arbitrary 
linear transformations of $V$, while the 
splitting of $\Lm{2}{V}$ into self-dual
and anti-self-dual 2-forms is not.
To see the result, note that a 2-plane
$P$ can be represented by a 2-vector
$u \wedge v$ with $u,v$ 
a basis for $P$, and this 
2-vector is unique up to scaling. With orientation
taken into account, it is unique up 
to positive rescaling. Conversely, every
2-vector $\gamma$ with $\gamma^2=0$ is
of this form for a unique 2-plane $P$.
Taking a metric on $V,$  we identify
2-vectors and 2-forms, and split our 2-vector
into self-dual and anti-self-dual pieces:
if $x^1,\dots,x^4$ are orthogonal coordinates on $V$,
and $dx^{ij} = dx^i \wedge dx^j,$ 
we can express a 2-form uniquely as
\begin{equation} \label{eqn:Klein} 
\begin{split}
\gamma & = \left ( X_1 + Y_1 \right ) dx^{12}
+  \left ( X_1 - Y_1 \right ) dx^{34}
+ \left ( X_2 + Y_2 \right ) dx^{31} \\
& + \left ( X_2 - Y_2 \right ) dx^{24}
+ \left ( X_3 + Y_3 \right ) dx^{23}
+  \left ( X^3 - Y_3 \right ) dx^{14}
\end{split}
\end{equation}
so that
\[
\gamma^2 = 2 \left ( X_1^2 + X_2^2 + X_3^3 - Y_1^2 - Y_2^2 - Y_3^2 \right)
dx^{1234}.
\]
The space of 2-forms satisfying $\gamma^2=0$ is just the
space of $X,Y$ with $|X|=|Y|$. To modulo out
by rescaling, we can always take $|X|^2+|Y|^2=2.$
Then it is clear that $\Gro{2}{V}=S^2_{+} \times S^2_{-}$
where $X$ provides coordinates on $S^2_{+}$ 
and $Y$ on $S^2_{-}$. We will now see that in these coordinates
the conformal quadratic form is 
$\left|\dot{Y}\right|^2-\left|\dot{X}\right|^2$.
Take the family of 2-planes $P(t)$ given
by
\[
\begin{pmatrix}
dx^3 \\
dx^4 
\end{pmatrix}
=-t \dot{p}
\begin{pmatrix}
dx^1 \\
dx^2
\end{pmatrix}
\]
and consider the 2-vector
\[
u(t) \wedge v(t) = 
\left( \partial_1 + t \dot{p}^1_1 \partial_3 + t \dot{p}^2_1 
\partial_4 \right ) \wedge 
\left ( \partial_2 + t \dot{p}^1_2 \partial_3 + t \dot{p}^2_2 
\partial_4 \right ) 
\]
which represents the 2-plane $P(t).$
Calculate that
\[
\left ( \left . \frac{d}{dt} u(t) \wedge v(t) \right |_{t=0} \right )^2
= - 2 \, \det \dot{p} \, dx^{1234}.
\]
\end{proof}

\begin{corollary} If $X \subset \Gro{2}{V}$ is
a compact elliptic line congruence, then
$X$ is diffeomorphic to $S^{2-}$,
hence a sphere. A compact elliptic line
congruence is precisely the graph of
a strictly contracting map
\[
S^{2-} \to S^{2+}
\]
\end{corollary}

For example, the Riemann sphere of a complex structure
$J : V \to V$ will be represented,
in complex linear coordinates on $V$,
by a constant map \( S^{2-} \to S^{2+}. \)
However, after a real linear transformation
of $V$, it can look somewhat different;
the picture of the Riemann sphere in
arbitrary real linear coordinates is
worked out in \cite{Gage:1985}.

\subsection{Real curves}
\begin{corollary}
Given an elliptic line congruence $X \subset \Gro{2}{V}$,
every line through the origin in $V$ is contained
in a discrete set of 2-planes belonging to $X$.
Similarly (replacing $V$ by $V^*$), 
every 3-plane in $V$ contains a discrete
set of 2-planes belonging to $X$.
\end{corollary}
\begin{proof}
The set of 2-planes containing the $x^1$ axis
is given in our coordinates by $X_1=Y_1,X_2=Y_2,X_3=-Y_3$
since these 2-planes must be represented
as 2-vectors by $\partial_1 \wedge (\dotsb)$.
But this is the graph of an isometry $S^2_{-} \to S^2_{+}$
and therefore cannot have any tangent vectors
in common with a strictly contracting map.
\end{proof}

The set of 2-planes containing a given
line sits inside $\Gro{2}{V}$ as the
graph of an orientation reversing
isometry $S^2_{-} \to S^2_{+},$ and
every orientation reversing isometry
occurs this way. To see this, use the
action of $\SO{4}$ on $\Gro{2}{V}.$
Similarly the orientation preserving
isometries are constructed by taking
the 2-planes contained in a given
3-plane.

The homology class of the set of 2-planes
containing a given line is therefore
$\left[S^2_{-}\right]-\left[S^2_{+}\right].$
The intersections with $S^2_{-} \times \text{pt}$
are therefore all negative, and so
are the intersections with any elliptic
line congruence because the intersections
do not change sign as we deform 
elliptic line congruences.

\begin{corollary}
A  compact elliptic
line congruence contains exactly one 2-plane
in $V$ containing each line through $0$ in $V$.
Moreover, the 2-plane depends smoothly on the
choice of line, by transversality.
\end{corollary}
\begin{proof} The intersections are all
negative, but the intersection number is $-1.$
\end{proof}

We think of this as infinitesimally solving
a Cauchy problem, complexifying real curves.
Replacing $V$ by $V^*$, we see that a compact
elliptic line congruence contains exactly one
2-plane in $V$ contained in each 3-plane
through $0$ in $V$. We think of this
as infinitesimal real hypersurface geometry.

\begin{corollary}
Any two 2-planes belonging to the same
compact elliptic line congruence are
transverse.
\end{corollary}

\subsection{Families of elliptic line congruences}
\begin{proposition} The graph of a 
strictly contracting map \( S^2 \to S^2 \)
has image contained in the interior of a single hemisphere.
\end{proposition}
\begin{proof} Suppose that there are two
antipodal points $x_1, x_2 \in S^2$ that get mapped to the
same point $y \in S^2$. Then any point in the image
must be the image of a point $x$ lying
in a hemisphere about either $x_1$ or $x_2$.
So any point in the image lies in the interior of 
a hemisphere about $y$. 
Therefore we need only show that some pair of antipodal
points $x_1, x_2$ get mapped to the same point $y$.

The image cannot contain
two antipodal points, because their preimages
could be no further apart than antipodal points,
and the map is strictly contracting. 
Therefore, the image misses some point
of $S^2$ and
the result is clear from the following lemma.
\end{proof}
\begin{lemma} Any continuous map
\[
\phi : S^2 \to \mathbb{R}^2
\]
must map some pair of antipodal points to
the same point.
\end{lemma}
\begin{proof} Suppose otherwise. Let
\[
\sigma(x) =
\frac{\phi(x)-\phi(-x)}{\left \| \phi(x)-\phi(-x) \right \|}.
\]
This is a continuous map
\(
\sigma : S^2 \to S^1
\)
so that 
\[
\sigma(-x)=-\sigma(x) \; .
\]
Because $S^2$ is simply connected, there is some
continuous map
\(
\tilde{\sigma} : S^2 \to \mathbb{R}
\)
so that
\[
\sigma = e^{i \tilde{\sigma} } \; .
\]
But then
\[
\tilde{\sigma}(-x)=\tilde{\sigma}(x)+(2k+1) \pi
\]
for some integer $k$. However, plugging
in $-x$:
\begin{align*}
\tilde{\sigma}(x) &= \tilde{\sigma}\left (- \left ( -x \right ) \right ) \\
              &= \tilde{\sigma}(-x)+(2k+1) \pi \\
              &= \tilde{\sigma}(x)+2(2k+1) \pi
\end{align*}
which is a contradiction.
\end{proof}

\begin{corollary} Every compact elliptic line congruence
can be smoothly deformed into any other through a family of
compact elliptic line congruences.
\end{corollary}
\begin{proof}
Rotation by $\SO{4}$ and dilation into a map
to a point takes us into the Cauchy--Riemann
equations.
\end{proof}

\begin{corollary} \label{cor:deformTame}
Every finite-dimensional
family $X \to B$
of elliptic line congruences in 
the fibers of a vector bundle $V \to B$ 
can be deformed globally and smoothly 
into a family of Riemann spheres, i.e.\
a complex vector bundle structure on $V.$
\end{corollary}
\begin{proof} Smoothly pick a positive definite
inner product on each fiber of $V.$ 
The ``center of mass''
of the image of each strictly
contracting map $S^2_{-} \to S^2_{+}$
provides a target to deform to.
\end{proof}

\subsection{Taming}

\begin{definition} A symplectic form \( \Omega \in \Lm{2}{V^*} \)
\emph{tames} an elliptic line congruence \( X \subset \Gro{2}{V} \)
if $\Omega$ is positive on every oriented 2-plane $P \subset V$
such that $P \in X$. An elliptic line congruence
is \emph{tameable} if it is tamed by some symplectic
form.
\end{definition}

\begin{corollary} \label{cor:Taming}
Every compact elliptic line congruence is tameable.
\end{corollary}
\begin{proof} Suppose that the line
congruence is written as the graph of a map
\[
\phi : S^{2-} \to S^{2+} \; .
\]
Suppose that the image of $\phi$ is contained
in the hemisphere about $\Omega \in S^{2+}$.
Suppose also that we have selected on $V$
not just a conformal class of positive definite
quadratic form, but an actual quadratic
form. Then our $\Omega$ can be thought of
as a linear function on $S^{2+}$ which is
positive on the image of $\phi$.
At the same time, $\Omega \in \Lm{2+}{V} \cong \Lm{2+}{V^*}$
is a 2-form. 
Take any $P \in X$, i.e.\ 
\[
P=(v_1 \wedge v_2, w_1 \wedge w_2) \in S^{2+} \times S^{2-}.
\]
We have 
\[
\Omega \in S^{2+} \subset \Lm{2+}{V^*} 
\]
so that
\[
\Omega(P) = \Omega(v_1 \wedge v_2) > 0
\]
because $\Omega$ is positive on the image of $\phi$.
Therefore $\Omega$ is a taming symplectic structure.
\end{proof}

\begin{lemma} If a 2-form is positive on every
2-plane belonging to a compact elliptic line congruence
(i.e.\ taming), then it is symplectic.
The space of symplectic forms
taming a compact elliptic line congruence
is an open convex cone in the space of symplectic
2-forms. 
\end{lemma}
\begin{proof}
If $\omega>0$ on every 2-plane in $X$ a compact
elliptic line congruence, then pick $P_1,P_2 \in X$
with $P_1 \ne P_2.$ Then $P_1 \cap P_2 = 0$
so that $V = P_1 \oplus P_2$ and $\omega>0$ on
both $P_1$ and $P_2$ so $\omega^2 > 0$ on $V$.

If two 2-forms $\omega_1,\omega_2$ both tame
a compact elliptic line congruence $X,$ then
clearly $a \omega_1 + b \omega_2$ does as
well for any $a,b>0.$
\end{proof}

\begin{definition} A pseudocomplex
structure $E \subset \Gro{2}{TM}$ is \emph{tamed}
by a symplectic form $\Omega \in \nForms{2}{M}$
if $\Omega_m \in \Lm{2}{T^*_m M}$ tames $E_m \subset \Gro{2}{T_m M}$
for every $m \in M$.
\end{definition}

\begin{corollary} Any proper pseudocomplex structure
is locally tameable.
\end{corollary}

The significance of taming symplectic forms
will soon be made clear.

\begin{proposition} The space of compact elliptic
line congruences tamed by a given symplectic
form \( \Omega \in \Lm{2}{V^*} \) is 
contractible.
\end{proposition}
\begin{proof} We can arrange, by linear
transformation, that $\Omega$ is
an element of $S^{2+}$. Then the taming
condition on an elliptic line congruence
is that, thinking of it as a map
\( \psi : S^{2-} \to S^{2+}, \)
its image is contained in the hemisphere
of $S^{2+}$ containing $\Omega$.
But then we can ``shrink'' \( \psi, \)
for example using the Riemannian
geometry exponential map about $\Omega$
in $S^{2+}$, which contracts the image
and preserves the ellipticity
(the strict contractivity) of $\psi$.
\end{proof}

A 2-form $\omega$ on a manifold $M^4$ 
with the property that $\omega^2$
does not vanish at any point is called
an \emph{almost symplectic structure}.

\begin{proposition} The space of almost
symplectic structures taming a fixed
proper pseudocomplex structure is an affine
space, not empty.
\end{proposition}
\begin{proof}
We can pick these almost symplectic structures
locally, and glue them together using convex
positive combinations.
\end{proof}

\begin{proposition} A four-dimensional manifold
admits a proper pseudocomplex structure  precisely
when it admits an almost complex structure,
which occurs precisely when it admits 
an almost symplectic structure.
\end{proposition}
\begin{proof}
As is well known, and proven in Steenrod
\cite{Steenrod:1999}, a manifold admits
an almost symplectic structure precisely when
it admits an almost complex structure,
and that almost complex structure
can be chosen to be tamed by the almost symplectic
structure. Then the Cauchy--Riemann 
equations of pseudoholomorphic curves
for the almost complex structure provide
a tamed proper pseudocomplex structure. 
\end{proof}

\begin{corollary} Every symplectic 4-manifold
admits a proper pseudocomplex structure
tamed by its symplectic structure.
\end{corollary}

\begin{corollary} Every proper pseudocomplex
structure tamed by a given almost symplectic
structure can be deformed through  proper
pseudocomplex structures tamed by that
same almost symplectic structure into an
almost complex structure.
\end{corollary}
\begin{proof} As for Corollaries 
\ref{cor:deformTame}
and \ref{cor:Taming}.
\end{proof}

\begin{proposition} Suppose that $X \subset \Gro{2}{V}$
is a compact elliptic line congruence.
Then any two distinct oriented 2-planes belonging to $X$ 
intersect transversely and positively 
at the origin.
\end{proposition}
\begin{proof} The positivity of 
intersection holds for Riemann spheres,
while transversality holds for all
compact elliptic line congruences,
and we can deform any compact
elliptic line congruence into a 
Riemann sphere.
\end{proof}

\begin{corollary} Any two transverse $E$ curves of
a proper pseudocomplex structure $E$
intersect positively.
\end{corollary}

A deeper problem is the positivity of
intersections of singular and nontransverse
$E$ curves. We will address this problem
in Section \ref{sec:GMY}.

\subsection{The Grassmannian as homogeneous space}
The Grassmannian is a homogeneous space 
\[
\Gro{2}{\R{4}} = \GLp{4,\R{}}/H_0
\]
where $H_0$ is the group of matrices of the
form
\[
\begin{pmatrix}
a & b \\
0 & d
\end{pmatrix}
\]
with $\det a, \det d >0$, i.e.\ the isotropy
subgroup of the 2-plane $P_0=\left(x^3=x^4=0\right).$
(Here, $\GLp{4,\R{}}$ means the $4 \times 4$ matrices
with positive determinant.) The tangent space
to $\Gro{2}{\R{4}}$ is identified with the
quotient of Lie algebras, i.e.\ with
matrices 
\[
\begin{pmatrix}
a & b \\
c & d
\end{pmatrix}
\]
modulo those of the form
\[
\begin{pmatrix}
a & b \\
0 & d
\end{pmatrix},
\]
hence with the $c$ component, i.e.\ the lower
left $2 \times 2$ matrix. To relate this
point of view to the others we have pursued,
consider a family of 2-planes $P(t)$
given by the equations
\[
\begin{pmatrix}
dx^3 \\
dx^4 
\end{pmatrix}
+ t \dot{p} 
\begin{pmatrix}
dx^1 \\
dx^2
\end{pmatrix}.
\]
Then $P(t)=g(t)P_0$ where
\[
g(t) = 
\begin{pmatrix}
1_2 & 0 \\
-t \dot{p} & 1_2
\end{pmatrix}.
\]
Here $1_2$ means the $2 \times 2$ identity
matrix. Thus in terms of the homogeneous
space point of view, $c = -\dot{p}$
and the conformal quadratic form is $\det  c = \det  \dot{p}.$

The action of $H_0$ on the tangent
space $T_{P_0} \Gro{2}{\R{4}}$ is given by
quotienting the adjoint action:
\[
\begin{pmatrix}
a & b \\
0 & d 
\end{pmatrix}
\begin{pmatrix}
* & * \\
c & * 
\end{pmatrix}
\begin{pmatrix}
a & b \\
0 & d
\end{pmatrix}^{-1}
=
\begin{pmatrix}
* & * \\
dca^{1} & *
\end{pmatrix}
\]
so the action is $c \mapsto dca^{-1}$
which clearly preserves the quadratic
form $c \mapsto \det c$ up to a positive
factor. Moreover the action is transitive
on positive tangent vectors to the 
Grassmannian, i.e.\ on matrices
$c$ with positive determinant, and on
negative vectors, and on nonzero null
vectors, i.e.\ matrices $c$ of rank 1.

Obviously the symmetric
bilinear form 
associated to the quadratic form $c \mapsto \det c$
is
\[
\left<
\begin{pmatrix}
A^1_1 & A^1_2 \\
A^2_1 & A^2_2
\end{pmatrix}
,
\begin{pmatrix}
B^1_1 & B^1_2 \\
B^2_1 & B^2_2
\end{pmatrix}
\right>
=
A^1_1 B^2_2 + B^1_1 A^2_2 - A^1_2 B^2_1 - B^1_2 A^2_1.
\]

\begin{lemma}
Consider bases of $T_{P_0} \Gro{2}{\R{4}}$ which
are orthonormal for the quadratic form,
up to a positive scaling factor, and are positively
oriented. The group $H_0$ acts on such bases with
two orbits.
\end{lemma}
\begin{proof}
Take such a basis $v_1,v_2,v_3,v_4.$ We may
rescale it to ensure that 
\[
\left<v_i,v_j\right>=
\begin{cases}
1 & \text{ if } i=j=1 \text{ or } 2 \\
-1 & \text{ if } i=j=3 \text { or } 4 \\
0 & \text{ otherwise.}
\end{cases}
\]
Think of each $v_i$ as a $2 \times 2$ matrix.
Using the action $v \mapsto dva^{-1}$
we can first arrange that $v_1=1_2$.
Now we have to work with the subgroup of 
$H_0$ fixing $1_2,$ which is the
group of transformations $v \mapsto ava^{-1}.$
Calculate that $\left<v_1,v\right>=\tr v$
so that $v_2,v_3,v_4$ are now traceless.
But then $v_2$ must be traceless, have
determinant 1, and so has minimal polynomial
$t^2+1$, i.e.\ is a complex structure.
We can therefore arrange by $v \mapsto ava^{-1}$
that $v_2= \pm J$ where
\[
J =
\begin{pmatrix}
0 & -1 \\
1 & 0
\end{pmatrix}.
\]
Let $K$ be complex conjugation, i.e.\
\[
K  =
\begin{pmatrix}
0 & 1 \\
1 & 0
\end{pmatrix}.
\]
Check that $K$ and $JK$ are perpendicular
to $v_1=1_2$ and $v_2=\pm J$, and to one another.
So $v_2$ and $v_3$ must be obtained
from $K$ and $JK$ by an orthogonal
transformation of the plane they span.
But we still have the freedom 
to employ the subgroup of $H_0$
fixing $v_1=1_2$ and $v_2=\pm J$, 
i.e.\ the transformations $v \mapsto ava^{-1}$
where $a$ is complex linear. This
enables us to rotate $v_3$ into $K$ and
then get $v_4 = \pm JK.$

But now the orientation of the basis
forces the two $\pm$ signs to be equal.
\end{proof}

\begin{corollary}
The general linear group $\GLp{V}$ acts
transitively on $\Gro{2}{V}$ and moreover
acts transitively 
on the positive definite 2-planes in the
tangent spaces of $\Gro{2}{V}.$
\end{corollary}
\begin{proof}
Given any positive definite 2-plane 
$\Pi \subset T_{P_0} \Gro{2}{V}$
we can take a basis $v_1,v_2,v_3,v_4$
for $T_{P_0} \Gro{2}{V}$ which
is orthonormal, up to a positive scaling
factor, with $v_1,v_2$ spanning $\Pi$ 
and we can use the group $H_0 \subset \GLp{V}$
to get it to the form $1_2,\pm J, K, \pm JK.$
\end{proof}

\begin{corollary}
A line congruence is tangent to
a Riemann sphere at a point precisely
when it is elliptic near that point.
\end{corollary}

\subsection{Totally real surfaces}

\begin{definition}
Take \( X \subset \Gro{2}{V} \)
an elliptic line congruence, and
$R \subset V$ a 2-plane,
so that with either orientation,
$R$ does not belong to $X$.
Call such an $R$ a \emph{totally real
$2$-plane} for $X$.
\end{definition}

\begin{lemma}
A vector $w$ belongs to a 2-plane $P$
represented by a 2-vector $u \wedge v$
precisely when $w \wedge u \wedge v = 0.$
\end{lemma}

\begin{lemma}
A pair of 2-planes $P_0$ and $P_1$
represented by 2-vectors $u_0 \wedge v_0$
and $u_1 \wedge v_1$ intersect in
a line or coincide (modulo orientation) precisely when
\[
u_0 \wedge v_0 \wedge u_1 \wedge v_1 = 0.
\]
\end{lemma}
\begin{proof}
If the 2-planes $P_0$ and $P_1$ contain
a common nonzero vector, we can take $u_0$ and $u_1$
to be that vector. Conversely, if 
$u_0 \wedge v_0 \wedge u_1 \wedge v_1 = 0,$
then by Cartan's lemma
\[
u_1 \wedge v_1 = u_0 \wedge x + v_0 \wedge y
\]
for some vectors $x$ and $y.$ 
Squaring both sides
\[
0 = u_0 \wedge v_0 \wedge x \wedge y.
\]
Without loss of generality, we take
$u_0=e_1,v_0=e_2$ in a basis $e_1,e_2,e_3,e_4.$
We find then that the matrix
\[
\begin{pmatrix}
x_3 & y_3 \\
x_4 & y_4
\end{pmatrix}
\]
has vanishing determinant, so has kernel.
Then take $v = v_1 e_1 + v_2 e_2$ and
\begin{align*}
v \wedge u_1 \wedge v_1 &= 
v \wedge \left ( e_1 \wedge x + e_2 \wedge y \right ) \\
&= e_1 \wedge e_2 \wedge \left ( v_1 y - v_2 x \right )
\end{align*}
so we can pick $v$ to make this vanish.
\end{proof}

Using equation \vref{eqn:Klein} to write 
$u_0 \wedge v_0$ in terms of $X_0$ and $Y_0,$
one finds
\[
u_0 \wedge v_0 \wedge u_1 \wedge v_1 
=
2 \left ( \left<Y_0,Y_1\right> - \left<X_0,X_1\right> \right )
dx^{1234}. 
\]
\begin{lemma}
A pair of 2-planes $P_0$ and $P_1$ represented
by points $\left(X_0,Y_0\right)$ and $\left(X_1,Y_1\right)$
in $S^2_{+} \times S^2_{-}$ intersect
in a line or coincide (modulo orientation)
precisely when
\[
\left< X_0, X_1 \right > = \left < Y_0, Y_1 \right >.
\]
\end{lemma}

\begin{proposition}
Let $X \to \Gro{2}{V}$ be an elliptic line
congruence and $R$ a totally real
2-plane for $X$.
Define \( X(R) \) to be the elements
of $X$ which, as
2-planes in $V$, intersect the 2-plane $R$ in
a line. Then \( X(R) \subset X \)
is a smooth curve (called the 
$R$ real points of $X$).
\end{proposition}
\begin{proof} 
Suppose that $X(R)$ has a singular
point. Pick a complex structure on 
$V$ tangent to $X$ at that singular
point.  Suppose that 
$R$ is represented in the splitting
\[
\Gro{2}{V} = S^{2+} \times S^{2-}
\]
by $R = (X_0, Y_0),$ while the line
congruence $X$ is 
represented, at least locally, by
a strictly contracting map $\psi : S^2_{-} \to S^2_{+}.$  
Then a point
$(\psi(Y),Y)$ on the graph of $\psi$ 
represents a 2-plane intersecting
the 2-plane $R$ in at least a line precisely when
\[
0 = - \left < X_0, \psi(Y) \right > 
+ \left < Y_0, Y \right >.
\]
We can arrange that the Riemann
sphere tangent to $X$ at this 
point be $X_0 \times S^2_{-}.$
Taking the gradient with respect to $x$, we
find that singular points will be precisely
those where
\[
Y_0 - \psi'(Y)^t X_0 \perp T_Y S^{2-}
\]
where $()^t$ indicates transpose. But by
the tangency condition, $\psi'(Y)=0$.
So therefore \( Y_0 = \pm Y. \) But then
\begin{eqnarray*}
0 & = & - \left < X_0, \psi(Y) \right > + \left < Y_0, Y \right >\\
  & = & - \left < X_0 , \psi(Y) \right > \pm 1 
\end{eqnarray*}
so that $X_0 = \pm \psi(Y)$. This shows
that $R$ is $(\psi(Y),Y)$ up to reorientation,
so that intersection occurs on more than
a line. Hence $R$ is not totally real.
\end{proof}

\begin{corollary} If $X \subset \Gro{2}{V}$ is
a compact elliptic line congruence and 
$R \subset V$ is a totally real 2-plane 
for $X$, then $X(R) \subset X$ is a smooth
embedded circle.
\end{corollary}
\begin{proof} We can deform $X$ to any
of its osculating complex structures
while keeping $R$ totally real. We may
have to move $R$ while we move $X,$
but it is easy to arrange that $R,$
with either orientation, never belongs
to $X$ during the deformation of $X$
since $R$ is just a point in the 
Grassmannian, while $X$ is a surface. In the
process, we never generate a singular
point on $X(R)$. The result is now a
calculation for the standard Riemann
sphere on $V = \mathbb{C}^2$.
\end{proof}

\subsection{The moving frame}
We will now employ Cartan's method of the
moving frame to uncover differential invariants
of elliptic line congruences, invariant
under orientation preserving linear transformations
of $V$.
Let $B_0$ be the set of orientation preserving
linear maps
identifying $V$ with $\R{4}$.
Pick a subspace $\R{2} = \R{2} \oplus 0  \subset \R{4}$.
Consider the map $B_0 \to \Gro{2}{V}$ 
given by taking for each $\lambda \in B_0$
the subspace $\lambda^{-1} \R{2} \subset V$,
and we  orient this subspace so that
$\lambda$ is orientation preserving.
Then $B_0 \to \Gro{2}{V}$ 
is a principal left $H_0$ bundle, where $H_0$ is
the group of invertible matrices of the form
\[
\begin{pmatrix}
a & b \\
0 & c
\end{pmatrix}
\]
with $a,b,c$ representing $2 \times 2$ matrices,
and $\det a > 0$ and $\det c > 0.$

  On $B_0$ we have the Maurer--Cartan 1-form
\[
\mu = d \lambda \, \lambda^{-1} \in \nForms{1}{B_0} 
\otimes \gl{4,\mathbb{R}}
\]
which satisfies
\begin{equation} \label{eqn:dmu}
d \mu = \mu \wedge \mu
\end{equation}
and under left action of $H_0$
\[
L_h^* \mu = \Ad_h \mu \; .
\]
We will split $\mu$ into complex linear
and complex conjugate linear parts:
\[
\mu \cdot v = 
\begin{pmatrix}
\xi & \eta \\
\vartheta & \zeta
\end{pmatrix}
v
+
\begin{pmatrix}
\xi' & \eta' \\
\vartheta' & \zeta'
\end{pmatrix}
\bar v
\]
for $v \in \R{4} = \C{2}$. Then we calculate
that equation (\ref{eqn:dmu}) becomes
\[
d 
\begin{pmatrix}
\xi \\
\xi' \\
\eta \\
\eta' \\
\vartheta \\
\vartheta' \\
\zeta \\
\zeta'
\end{pmatrix}
=
\begin{pmatrix}
\xi' \wedge \bar \xi' + \eta \wedge \vartheta 
+ \eta' \wedge \bar \vartheta' \\
\xi \wedge \xi' + \xi' \wedge \bar \xi 
+ \eta \wedge \vartheta' + \eta' \wedge \bar \vartheta \\
\xi \wedge \eta + \xi' \wedge \bar \eta' + \eta \wedge \zeta 
+ \eta' \wedge \bar \zeta' \\
\xi \wedge \eta' + \xi' \wedge \bar \eta + \eta \wedge \zeta'
+ \eta' \wedge \bar \zeta \\
\left ( \zeta - \xi \right ) \wedge \vartheta 
- \bar \xi' \wedge \vartheta' + \zeta' \wedge \bar \vartheta' \\
- \xi' \wedge \vartheta + \zeta' \wedge \bar \vartheta
- \left ( \bar \xi - \zeta \right ) \wedge \vartheta' \\
\vartheta \wedge \eta + \vartheta' \wedge \bar \eta' 
+ \zeta' \wedge \bar \zeta' \\
\vartheta \wedge \eta' + \vartheta' \wedge \bar \eta 
+ \zeta \wedge \zeta' + \zeta' \wedge \bar \zeta
\end{pmatrix} \; .
\]
On each fiber of $B_0 \to \Gro{2}{V}$, we have
$\vartheta = \vartheta' = 0$ , i.e.\ $\vartheta, \vartheta'$
are semibasic.

Now suppose that we have an elliptic
line congruence $X \subset \Gro{2}{V}$.
First, we can form the pullback bundle
\[
\xymatrix{
B_1(X) \ar[r] \ar[d] & B_0 \ar[d] \\
X \ar[r] & \Gro{2}{V}.
}
\]
This is also a principal left $H_0$ subbundle.
The fibers of $B_0 \to \Gro{2}{V}$ are
cut out by the equations 
\[
\vartheta = \bar \vartheta  = \vartheta' = \bar \vartheta' = 0
\]
(four independent equations)
so that $\vartheta, \bar \vartheta, \vartheta', \bar \vartheta'$
are semibasic for $B_0 \to \Gro{2}{V}$. 
The fibers of $B_1(X)$ are cut out by the
same equations, but on $B_1(X)$. Since $X$ has only 
two dimensions,
we must have at least two relations among
these four 1-forms on $B_1(X)$. 
Let us try to normalize these relations.

To do this, let us examine $X$ near some point
$p_0 \in X$.
This surface $X$ is tangent to a Riemann
sphere at each point, so pick 
a Riemann sphere $X_0$ 
tangent to $X$ at $p_0$.
Suppose that $X_0$ is the
Riemann sphere of the complex structure
\(
J_0 : V \to V \; .
\) 
Let $H_2$ be the subgroup of $H_0$
consisting of matrices
\[
\begin{pmatrix}
a & b \\
0 & c 
\end{pmatrix}
\] 
where $a,c \in \mathbb{C}^{\times}$ are
complex linear, and $b$ is still
a real linear $2 \times 2$ matrix.
The tangent space 
\( T_{p_0} X_0  \) to a Riemann sphere
associated to a complex structure 
$J_0 : V \to V$
is canonically identified with the
$J_0$ complex linear maps $p_0 \to V/p_0$.
Above each point $p$
in our Riemann sphere, we have
a distinguished subset of 
$B_1 \left ( X_0 \right )$:
those elements 
$\lambda \in B_1\left ( X_0 \right )$ which 
identify $p = \lambda^{-1} \R{2}$
with $\R{2} = \mathbb{C}$
by a complex linear map.
Call this $B_2 \left ( X_0 \right )$.
It is easy to show by homogeneity under 
$\GL{2,C}$ that $B_2 \left ( X_0 \right )$
is a manifold and that 
$B_2 \left ( X_0 \right ) \to X_0$
is a principal left $H_2$  bundle.
If $X_0$ and $X_1$ are
the Riemann spheres of complex
structures $J_0$ and $J_1$, and
$X_0$ is tangent to $X_1$ 
at a point $p \in \Gro{2}{V}$, then
we have canonical identifications:
\[
\left \{ J_0 \text{ linear maps } p \to V/p \right \}
\cong
T_p X_0
=
T_p X_1
\cong
\left \{ J_1 \text{ linear maps } p \to V/p \right \} \; .
\]
\begin{lemma} If $X_0$ and $X_1$ are
the Riemann spheres of complex structures
$J_0$ and $J_1$ on a four-dimensional vector
space $V$, and $X_0$ is tangent to $X_1$
at a point $p \in \Gro{2}{V}$, then 
$J_0$ and $J_1$ agree on $p$ and on $V/p$.
\end{lemma}
\begin{proof}
The characteristic variety of the tableaux
\[
A_k = \left \{ J_k \text{ linear maps } p \to V/p \right \}
\]
for $k=0,1$ is the same for $k=0$ as for $k=1$
because $A_0 = A_1$. But by definition,
\[
\Xi_{\C{}} \left ( A_k \right )
=
\left \{
v \in p \otimes \C{} \, | \,
\exists \xi \in \left ( (V/p) \otimes \C{} \right )^*,
\xi \otimes v \in A_k
\right \}.
\]
Calculating in complex coordinates,
we find that the characteristic
variety is the union of the two eigenspaces of $J_0$,
with eigenvalues $\pm \sqrt{-1}$.
This splitting of the characteristic
variety (as a real algebraic variety)
determines $J_0$ up to sign on $p$, and
so $J_1 = \pm J_0$ on $p$. If we had
a minus sign here, that would give
$p$ the opposite orientation. Similarly,
by taking transposes, we find
$J_1 = J_0$ on $V/p$. 
\end{proof}

This shows that the fibers match:
\[
B_2 \left ( X_0 \right )_p = B_2 \left ( X_1 \right )_p.
\]
We can now unambiguously
define the set
\[
B_2 \left ( X \right )
\]
for any elliptic line congruence $X$
by the equation
\[
B_2 \left ( X \right )_p
=
B_2 \left ( X_0 \right )_p
\]
whenever $X_0$ is a Riemann sphere
tangent to $X$ at $p$. It is
still not clear that this set 
$B_2 \left ( X \right )$ is a
principal left $H_2$ subbundle of $B_1(X) \to X.$

\begin{lemma}
$B_2 \left ( X \right ) \subset B_0$
is a smooth submanifold satisfying the equation
$\vartheta'=0$. 
\end{lemma}
\begin{proof} Recall that we defined
$B_0$ to be the real linear isomorphisms
\[
V \to \R{4}
\]
matching orientation. There is a map
\[
\lambda : B_0 \to \Lin{V}{\R{4}}
\]
defined by inclusion. The Maurer--Cartan
1-form is just $\mu = d \lambda \, \lambda^{-1}$.
Take $X_0$ a Riemann sphere tangent
to $X$ at $p$, and $J_0$ the associated
complex structure.
Let us assume that $V = \R{4} = \C{2}$, 
with complex linear coordinates $z,w$ for simplicity
of notation, so that $J_0$ is the
standard complex structure,
and consider the quantity
\[
Q = 
i \lambda^{-1} \left ( \partial_z \wedge \partial_{\bar z} \right )
\]
so that
\[
Q : B_0 \to \Lm{1,1}{V^*}.
\]
We calculate that
\[
J_0 Q = Q \text{ at } \lambda = I
\]
and that again at $\lambda = I$
\[
d (J_0 Q - Q) = 2i \vartheta' \partial_z \wedge \partial_w
+ 2i \bar \vartheta' \partial_{\bar w} \wedge \partial_{\bar z}
\]
so that indeed $\vartheta'=0$ precisely along the
directions of $T_{\lambda} B_2 \left ( X_0 \right )$,
where $X_0$ is the Riemann sphere of
our complex structure.

  Consequently, every Riemann sphere
$X_0$ tangent to $X$ at a point
$\lambda$ must satisfy
\[
T_{\lambda} B_2 \left ( X_0 \right ) =
\{ \vartheta' = 0 \} \subset T_{\lambda} B_1 \left ( X_0 \right )
=
T_{\lambda} B_1 \left ( X \right ) .
\]
This shows us that the relations among
$\vartheta, \bar \vartheta, \vartheta', \bar \vartheta'$
near the locus $B_2 \left ( X \right ) \subset
B_0 \left ( X \right )$ cannot
be
\[
\vartheta = \bar \vartheta = 0
\]
since these relations would also have to
hold on 
\[
T_{\lambda} B_2 \left ( X_0 \right )
=
T_{\lambda} B_2 \left ( X \right ).
\]
Indeed, on $B_2 \left ( X_0 \right )$
we still have $\vartheta \wedge \bar \vartheta \ne 0$,
since $\vartheta' = \bar \vartheta' = 0$ there,
but
\[
B_2 \left ( X_0 \right ) \to X_0
\]
is a submersion, and $\vartheta \wedge \bar \vartheta$
is semibasic for it, providing a basis for
the semibasic forms (check this in complex
coordinates). Therefore, near $B_2 \left ( X \right )$
we must find that the relations among 
$\vartheta, \bar \vartheta, \vartheta', \bar \vartheta'$
can be written
\[
\begin{pmatrix}
\vartheta' \\
\bar \vartheta' \\
\end{pmatrix}
=
\begin{pmatrix}
a & b \\
\bar b & \bar a
\end{pmatrix}
\begin{pmatrix}
\vartheta \\
\bar \vartheta
\end{pmatrix}.
\]
We can follow how these relations behave
as we travel up the fibers of $B_1 \left ( X \right )$.
We find that the coefficients $a,b$ can be
made to vanish, on a principal left $H_2$ bundle.
Moreover, it is easy to see that $B_2 \left ( X \right )$
is precisely this bundle, following the general
formalism presented in \cite{Jensen:1977}.
\end{proof}

Now we find that on $B_2 \left ( X \right )$, 
our equations are \( 0 = \vartheta' \) which implies
from \( d \mu = \mu \wedge \mu \) that
\[
\begin{pmatrix}
\xi' \\
\zeta' 
\end{pmatrix}
=
\begin{pmatrix}
f & h \\
-h & g 
\end{pmatrix}
\begin{pmatrix}
\vartheta \\
\bar \vartheta
\end{pmatrix}
\]
for some uniquely determined complex
valued functions
\[
f,g,h  : B_1 \to \mathbb{C} \; .
\]
Under the action of $H_2$ we find
that traveling up the fibers of 
$B_2 \left ( X \right ) \to X$,
\[
\begin{pmatrix}
a & b \\
0 & c
\end{pmatrix} \cdot 
\begin{pmatrix}
f \\
g \\
h
\end{pmatrix}
=
\begin{pmatrix}
a^2 \bar a^{-1} c^{-1} f \\
\bar a c \bar c^{-2} g \\
\left ( ah + b' \right ) \bar c^{-1}
\end{pmatrix}
\]
where $b'$ is the complex
conjugate linear part of $b$.
Therefore, letting $H_3$ be the group
of matrices
\[
\begin{pmatrix}
a & b \\
0 & c
\end{pmatrix}
\]
with $a,b,c$ complex numbers,
and $a, c \ne 0$,
we find a principal
left $H_3$ subbundle of $B_2 \left ( X \right )$,
call it $B_3 \left ( X \right )$, on which
$h=0$.

\begin{theorem} Let $V$ be a four-dimensional
real vector space, and $B_0$
the set of all linear isomorphisms $V \to \R{4}$.
Let $B_0 \to \Gro{2}{V}$ be the map
\[
\lambda \mapsto \lambda^{-1} \R{2} \oplus 0 .
\]
Let $H$ be the group of complex $2 \times 2$
matrices of the form
\[
\begin{pmatrix}
a & b \\
0 & c
\end{pmatrix}.
\]
Each elliptic line congruence 
\[
X^2 \subset \Gro{2}{V}
\]
determines invariantly a principal
left $H$ subbundle $B \left ( X \right ) \subset B_0$ so
that
\[
\xymatrix{
B (X) \ar[r] \ar[d] & B_0 \ar[d] \\
X \ar[r] & \Gro{2}{V}
}
\]
and so that the Maurer--Cartan 1-form 
$\mu = d \lambda \, \lambda^{-1}$
on $B_0$, when written out as
\[
\mu \cdot v = 
\begin{pmatrix}
\xi & \eta \\
\vartheta & \zeta
\end{pmatrix}
v
+
\begin{pmatrix}
\xi' & \eta' \\
\vartheta' & \zeta'
\end{pmatrix}
\bar v
\]
for $v \in \R{4} = \C{2}$,
satisfies on $B \left ( X \right )$
the equations
\begin{align*}
\vartheta' &= 0 \\
\xi' &= f \vartheta \\
\zeta' &= g \bar \vartheta \\
\eta' &= - s \vartheta - t \bar \vartheta
\end{align*}
for uniquely determined complex valued functions 
\[
f,g,s,t : B \left ( X \right ) \to \C{},
\]
and so that the covariant derivatives
\begin{align*}
\nabla f &= df - f \left ( 2 \xi - \bar \xi - \zeta \right ) \\
\nabla g &= dg - g \left ( \zeta - 2 \bar \zeta + \bar \xi \right )
\end{align*}
satisfy
\[
\nabla
\begin{pmatrix}
f \\
g
\end{pmatrix}
=
\begin{pmatrix}
u & s \\
t & v 
\end{pmatrix}
\begin{pmatrix}
\vartheta \\
\bar \vartheta
\end{pmatrix}
\]
for uniquely determined complex valued functions
\[
u, v : B(X) \to \C{} \; .
\]
\end{theorem}

At each point $p$ of $X$, we therefore
have a well-defined family of maps $\lambda : V \to \R{4}$,
forming the fiber of $B \left ( X  \right )$,
determined up to $H$ action. By complex 
linearity of the elements of the group $H$,
this determines a complex structure $J(p)$ on $V$,
for which the real 2-plane $p$ is a complex
line, and for which $V/p$ is also a complex
line. We call this the \emph{osculating complex
structure} to $X$ at $p$. It determines
a map
\[
X \to \Cstrucs{V} \; .
\]

The osculating complex structure
determines a complex structure on
$X$ itself in the obvious manner:
each tangent space to $X$ is
identified with the complex linear maps
\[
p \to V/p
\]
for the osculating complex structure,
hence a one-dimensional complex vector space.
Consider the map $X \to \Cstrucs{V}$ from our
surface $X$ to the space of complex
structures on the vector space $V$.
Recall that the inclusion 
$\Cstrucs{V} \subset \operatorname{sl}(V)$
is a coadjoint orbit of $\SL{V}$ under
the identification
\[
\operatorname{sl}(V) \cong \operatorname{sl}(V)^*
\]
given by the Killing form. Thus
$\Cstrucs{V}$ is a $\SL{V}$ homogeneous space,
preserving the symplectic structure given
by the Kirillov symplectic
structure on a coadjoint orbit, and
preserving the pseudo-K\"ahler structure
given by the Killing form. This provides
invariants of $X$, given by pulling back
these structures. Using the complex
structures on $X$ and $\Cstrucs{V}$, 
and the map $X \to \Cstrucs{V}$,
we can take $\partial$ and $\bar \partial$
of this map as invariants of an
elliptic line congruence. Essentially,
these are $f$ and $g$. It turns
out that if $X$ is compact, then 
$X \to \Cstrucs{V}$ is
never holomorphic, unless it is
constant, so that $X$ is the Riemann
sphere of a single complex structure
$J : V \to V$.  The symplectic form
on $\Cstrucs{V}$ pulls back to
\[
i \left ( |f|^2 - |g|^2 \right ) \vartheta \wedge \bar \vartheta  
\]
which consequently integrates to zero, giving
\begin{equation}
\int i |f|^2 \vartheta \wedge \bar \vartheta
=
\int i |g|^2 \vartheta \wedge \bar \vartheta \label{eqn:Blue}
\end{equation}
if $X$ is compact,
so that, roughly speaking, the map
\[
X \to \Cstrucs{V}
\]
is equal parts holomorphic and conjugate holomorphic.
Vanishing of the $f$ and $g$ invariants
of an elliptic line congruence occurs precisely
for Riemann spheres of complex structures.

\begin{proposition}
There are no compact homogeneous
elliptic line congruences except Riemann
spheres.
\end{proposition}
\begin{proof}
Let $X$ be a compact homogeneous
elliptic line congruence.
Then by homogeneity, the 2-forms $i f \vartheta \wedge \bar \vartheta$ 
and $i g \vartheta \wedge \bar \vartheta$ (which
one can calculate are invariantly defined on $X$ itself) must
be either everywhere vanishing or everywhere
not vanishing. Moreover, if one vanishes everywhere,
then by equation (\ref{eqn:Blue}) they both do.
But if neither vanishes anywhere, then
there is a distinguished subbundle of $B(X)$
on which $f=1$. Using the Maurer--Cartan
equations, we see that this renders $\vartheta$
well defined on $X$, and thereby forces 
$X$ to have a globally defined nowhere
vanishing 1-form, in fact a parallelism,
so that $X$ cannot be topologically a sphere.
\end{proof}

\subsection{Elliptic line congruences as nonlinear
complex structures}
\begin{proposition}
Given \( X \subset \Gro{2}{V} \)
a compact elliptic line congruence,
we can define a map
\( 
J_X : V \backslash 0 \to V \backslash 0
\)
by asking that on each 2-plane $P \subset V$
belonging to $X$, $J_X$ acts as the
osculating complex structure on $P$.
This defines a smooth map, which satisfies
\begin{align}
&J_X^2 = -1 \tag{$1$} \\
&J_X (tv) = t J_X (v) \tag{$2$}
\intertext{
for \( v \in   V \backslash 0 \) and 
\( t \in \R{} \backslash 0 \)
and}
&J_X \text{ restricts to each 2-plane } 
\Span{v,J_X v} 
\text{ to be a linear map. } \tag{$3$}
\end{align}

Conversely, given such a map, define
$X$ to be the set of 2-planes invariant
under the map, and orient these 2-planes
so that $v \wedge J v$ is positive. Then
$X$ is an elliptic line congruence.
\end{proposition}
\begin{proof}
We only have to prove that starting from some
such $J$
we can define $X$ as above, and it turns
out to be an elliptic line congruence.
It is clear that $X$ is compact, since
it is a subset of a Grassmannian (which
is compact) and defined by a closed
condition. Define a vector field
$Z$ on the sphere $V/\R{+}$ by
\[
Z(v) = Jv \mod \Span{v}.
\] 
This is a nowhere vanishing
vector field defining a circle action
on the sphere. Moreover, on the
2-plane $P$ containing a vector $v$
and the vector $Jv,$ this vector
field is just $Z(v)=J_P v$, where
$J_P$ is the osculating complex structure.
So the flow of $Z(v)$ on $P$ is
$t \mapsto e^{t J_P} v,$ and in
particular the flow curves on the
sphere $V/\R{+}$
are great circles. The quotient space
is clearly $X$. By properness and
freedom of the action, $X$ is
a smooth compact surface.

We have to show that $X$ is elliptic.
This follows from the result of 
Gluck and Warner \cite{Gluck/Warner:1983}
that great circle fibrations
of the 3-sphere correspond precisely
to elliptic line congruences. We won't
need this result, so the reader will
forgive the author for not providing
a complete proof of the results of Gluck and
Warner.
\end{proof}

\section{Applying Cartan's method of equivalence}

So far, we have only studied the microlocal
geometry. The full geometry of a pseudocomplex
structure
$E^6 \subset \Gro{2}{TM}$ requires Cartan's
equivalence method. The algorithm for this
method is explained thoroughly in \cite{Gardner:1989};
the calculations for this specific equivalence
problem are given in great detail
in my thesis \cite{McKay:1999}; therefore
only the result will be presented here. First,
to define a $G$ structure associated to
a pseudocomplex structure
\( E^6 \subset \Gro{2}{TM} \) on a 4-manifold 
$M$, we work on $E$ rather than working on $M$. Already
we know that the fibers
of $E \to M$ are elliptic line congruences,
so they have complex structures, and each
point of $E_m$ imposes a complex structure
on $T_m M$, the osculating complex structure.
Also, $E$ is equipped with a 4-plane field:
if we write $p : E \to M$ for the projection
\( E \subset \Gro{2}{TM} \to M \), then
we have 4-planes: 
\[
\Theta_e = p'(e)^{-1}(e) \subset T_e E
\]
where $e \in E$ is thought of as 
a 2-plane $e \subset T_m M$. Contained
in these 4-planes are the 2-planes
\[
V_e = \operatorname{ker} p'(e) \subset \Theta_e 
\]
which are the tangent spaces to the fibers of 
$E \to M$, so are tangent spaces of elliptic
line congruences, and hence have complex structure.
But also 
\[
\Theta_e / V_e \cong T_m M
\]
has the osculating complex structure.

  Thus at each point of $E$, we have
a 2-plane $V_e$ contained in a 4-plane $\Theta_e$,
and complex structures on $V_e$ and on $\Theta_e/V_e$.
We let $W_0$ be any six-dimensional vector space,
equipped with a 4 dimensional subspace $\Theta_0 \subset W_0$
and a 2 dimensional subspace $V_0 \subset \Theta_0$,
with complex structures on $V_0$ and $\Theta_0/V_0$,
and let $G$ be the group of linear transformations
$W \to W$ preserving this structure. Then
any pseudocomplex structure 
\( E^6 \subset \Gro{2}{TM} \) is canonically
equipped with a $G$ structure. Moreover, this
$G$ structure encodes completely the
pseudocomplex structure.

  This $G$ structure is input to Cartan's
method of equivalence, and after a single
step (without making any prolongations
of structure equations) we obtain an
$H$ structure, where $H \subset G$
is the group of matrices preserving
$V_0 \subset \Theta_0 \subset W$ and
preserving a complex structure on
$W$ for which $\Theta_0$ and $V_0$ are
complex subspaces, and preserving
a certain 2-form from \( \Lm{2}{\Theta_0} \otimes W / \Theta_0 \)
so that $H$, in complex coordinates $z_1,z_2,z_3$, is the
collection of matrices:
\[
\begin{pmatrix}
a & 0 & 0 \\
b & c & 0 \\
d & e & ac^{-1}
\end{pmatrix}
\]
with $V_0 = \{ z_1 = z_2 = 0 \}$ and $\Theta_0 = \{ z_1 = 0 \}$
and the 2-form is \( dz_2 \wedge dz_3 \mod \Theta_0 \).
Write the total space of this $H$ structure as
$B \to E$, so that $B$ is a collection of linear
maps identifying tangent spaces of $E$ with $W$,
and is a principal right $H$ bundle.
The Cartan structure equations are
\begin{equation} \label{eqn:struc}
d \begin{pmatrix}
\theta \\
\omega \\
\pi
\end{pmatrix}
=
-
\begin{pmatrix}
\alpha & 0 & 0 \\
\beta & \gamma & 0 \\
\delta & \varepsilon & \alpha - \gamma
\end{pmatrix}
\wedge
\begin{pmatrix}
\theta \\
\omega \\
\pi
\end{pmatrix}
-
\pi \wedge 
\begin{pmatrix}
\omega \\
\sigma \\
0
\end{pmatrix}
+
\begin{pmatrix}
\tau_1 \\
\tau_2 \\
\tau_3
\end{pmatrix}
\wedge \bar \theta
\end{equation}
where $\theta, \omega, \pi, \tau_1, \tau_2, \tau_3, \sigma$ 
are complex valued
1-forms on $B$, uniquely determined, and 
$\alpha, \beta, \gamma, \delta, \varepsilon$ 
are complex valued 1-forms, not uniquely
determined, and
\begin{align*}
0 &= \sigma \wedge \bar \theta \wedge \bar \omega \\
0 &= \tau_j \wedge \bar \omega \wedge \bar \pi
\end{align*}
and \( \theta \wedge \bar \theta \wedge \omega \wedge
\bar \omega \wedge \pi \wedge \bar \pi \ne 0 \). 
Consequently, we can write
\[
\begin{pmatrix}
\tau_1 \\
\tau_2 \\
\tau_3 \\
\end{pmatrix}
=
\begin{pmatrix}
T_2 & T_3 \\
U_2 & U_3 \\
V_2 & V_3
\end{pmatrix}
\begin{pmatrix}
\bar \omega \\
\bar \pi
\end{pmatrix}
\]
and
\[
\sigma = S_1 \bar \theta + S_2 \bar \omega
\]
and these $S,T,U,V$ are the lowest
order invariants of pseudocomplex structures. 

  Let $(\theta_0, \omega_0, \pi_0)$ 
be a local section of the $H$ bundle, so that
\[
\theta_0, \omega_0, \pi_0 \in \nForms{1}{E} \otimes \mathbb{C} \; .
\]
Then the 4-plane field $\Theta$ is cut out by
the equation $\theta_0=\bar \theta_0=0$,
while $V$ is $\theta_0=\bar \theta_0=\omega_0=\bar \omega_0=0$. 
The $H$ structure preserves an almost complex
structure $J$ on $E$, for which $\theta_0,\omega_0,\pi_0$
are $(1,0)$ forms. The 4-plane $\Theta$ is a field
of $J$-complex 2-planes, while $V$ is a field
of $J$-complex lines.
Every $E$ curve
$\Sigma \subset M$ has tangent planes belonging
to $E$, forming a surface called the \emph{prolongation}
of $\Sigma$. 

\begin{proposition} The prolongation $\hat \Sigma \subset E$
of an $E$
curve $\Sigma \subset M$ is a $J$-holomorphic
curve $\hat \Sigma \subset E$ tangent to the 4-plane 
field $\Theta$. Conversely, a surface
$\hat \Sigma \subset E$ which is tangent to the 4-plane
field $\Theta$ and transverse to the
fibers of $E \to M$ is locally the prolongation 
of an $E$ curve $\Sigma \subset M$.
\end{proposition}
\begin{proof} First, take $\Sigma \subset M$ 
an $E$ curve. The prolongation $\hat \Sigma \subset E$
is the set of tangent planes of $\Sigma \subset M$.
We find that $p : E \to M$ restricts to $p : \hat \Sigma \to \Sigma$.
So $p'(e) : T_e \hat \Sigma \to T_p \Sigma$ and thus
\[
T_e \hat \Sigma \subset p'(e)^{-1} T_p \Sigma. 
\]
Therefore if $x \in \Sigma$ and $e \in \hat \Sigma$ 
with $p(e)=x$, then $e = T_p \Sigma$. We have
\[
\Theta_e = p'(e)^{-1}e = p'(e)^{-1} T_p \Sigma 
\]
so
\( T_e \hat \Sigma \subset \Theta_e. \)
We have to show that $\hat \Sigma$ is $J$-holomorphic.
Choose (locally) a coframing
\[
\theta_0, \omega_0, \pi_0 \in \nForms{(1,0)}{E}
\]
from our $H$ structure.
Since $\Theta$ is $\theta_0=\bar \theta_0=0$, 
we have $\theta_0=\bar \theta_0 = 0$
on $\hat \Sigma$, and so $0 = d \theta_0 = d \bar \theta_0$.
But from the structure equations, this gives
\[
0 = \pi_0 \wedge \omega_0 = \bar \pi_0 \wedge \bar \omega_0 \; .
\]
Therefore the $(1,0)$ forms $\theta_0,\omega_0,\pi_0$
restrict to $\hat \Sigma$ to satisfy some
complex linear relations, and thus $\hat \Sigma$
is $J$-pseudoholomorphic.

Now take $\hat \Sigma \subset E$ any $\Theta$ tangent
surface transverse to the
fibers of $p : E \to M$. Since our desired result is local,
we can suppose that
\[
p : \hat \Sigma  \to \Sigma = p \left ( \hat \Sigma \right )
\]
is a diffeomorphism. Let $e \in \hat \Sigma, m=p(e) \in \Sigma$.
By definition,
\[
T_e \hat \Sigma \subset \Theta_e = p'(e)^{-1} e
\]
so that
\[
p'(e) T_e \hat \Sigma \subset e \subset T_m M
\]
and since $p : \hat \Sigma \to \Sigma$ is a diffeomorphism:
\[
p'(e) T_e \hat \Sigma = T_p \Sigma \; .
\]
Therefore, $e = T_p \Sigma$, and so $\hat \Sigma$
is the prolongation of $\Sigma$.
\end{proof}

\begin{proposition} 
\label{prp:Rebuild}
Suppose that we have some
1-forms $\theta, \omega, \pi, \alpha, \ldots$
and functions on a six-dimensional manifold satisfying the 
structure equations \vref{eqn:struc}, and that
\[
\theta \wedge \bar \theta 
\wedge \omega \wedge \bar \omega 
\wedge \pi \wedge \bar \pi \ne 0 \; .
\]
Then locally these define a pseudocomplex
structure on a 4-manifold.
\end{proposition}
\begin{proof} The 4-manifold is constructed
by taking the foliation
\[
\theta = \bar \theta = \omega = \bar \omega  = 0 
\]
which is locally a fiber bundle, and
producing a base space for that fiber
bundle. The rest is elementary, following
the general pattern of the equivalence method.
\end{proof}

We can prolong the exterior differential system
on $E$ given by the 4-plane $\Theta$, to form
$E^{(1)}$, and prolong that to $E^{(2)}$, and so
on. Each of these $E^{(k)}$ has an almost complex
structure, easily seen from the structure
equations on the $E^{(k)}$, which can be easily 
calculated (although the calculation 
is somewhat lengthy). The prolongations of 
an $E$ curve to all orders will be pseudoholomorphic
curves in each of the $E^{(k)}$.

\subsection{Microlocal invariants}

\begin{proposition} Take $E$ a pseudocomplex structure
on a 4-manifold $M$, and let $B \to E$ be the
induced $G$ structure constructed above. Call
this bundle $B$ the pseudocomplex structure
bundle of $E$. Now
take a point $m \in M$ and look at the elliptic
line congruence $E_m \subset \Gro{2}{T_mM}$.
We are faced with two principal bundles over 
$E_m$: (1)
\( B \left ( E_m \right ) \to E_m \) the 
bundle associated
to the elliptic line congruence, following
our results on line congruences above (call
it the line congruence bundle of $E_m$), and (2) the
bundle $B_m \to E_m$ given as the part of
the pseudocomplex bundle
that lies above $E_m$. Given a coframe 
\[
\left ( \theta_0, \omega_0, \pi_0 \in \Lm{1}{T_e E} \right )
\]
from the pseudocomplex bundle, we can map
it to a linear map
\[
u : T_m M \to \R{4}
\] 
belonging to the line congruence bundle
as follows: since the fiber of $E \to M$
is given by the equations
\[
\theta_0 = \bar \theta_0 = \omega_0 = \bar \omega_0 = 0
\]
we can treat $\theta_0, \omega_0$ as well-defined 
complex valued 1-forms on $T_m M$. Then we can
use the map
\[
u(w) = 
\begin{pmatrix}
\omega_0 \\
\theta_0
\end{pmatrix} \; .
\]
This gives an equivariant bundle map
\[
\xymatrix{
B_m \ar[rr] \ar[dr] & & B \left ( E_m \right  ) \ar[dl] \\
& E_m & 
}
\]
over each line congruence $E_m$. Under this map
\[
\begin{pmatrix}
\xi & \eta \\
\vartheta & \zeta
\end{pmatrix}
=
\begin{pmatrix}
- \gamma & - \beta \\
- \pi & - \alpha
\end{pmatrix}
\qquad
\begin{pmatrix}
\xi' & \eta' \\
\vartheta' & \zeta'
\end{pmatrix}
=
\begin{pmatrix}
- S_2 \pi & U_3 \bar \pi - S_1 \pi \\
0 & T_3 \bar \pi
\end{pmatrix}.
\]
Moreover
\[
\begin{pmatrix}
f \\
g 
\end{pmatrix}
= 
\begin{pmatrix}
-S_2 \\
T_3 
\end{pmatrix}
\]
\end{proposition}

\subsection{Identifying almost complex structures}

From our identification of the microlocal
invariants we have:
\begin{proposition}
The following are equivalent:
\begin{enumerate}
\item \( E \subset \Gro{2}{TM} \) is an almost complex
structure. 
\item
\( 0 = \sigma \wedge \bar \theta = \tau_1 \wedge \bar \omega \)
\item
\( 
0 = \sigma  = \tau_1 \wedge \bar \omega = \tau_3 \wedge \bar \omega
\)
\item The projection $E \to M$
is holomorphic for some almost complex structure
on $M$ (which is then necessarily given by $E$ itself,
as a subset of $\Gro{2}{TM}$).
\end{enumerate}
\end{proposition}

\begin{proposition}
\label{prop:WhenAC}
A proper pseudocomplex
structure is an almost complex structure
precisely when either
\( 0 = \sigma \wedge \bar \theta \)
or 
\( 0 = \tau_1 \wedge \bar \omega. \)
\end{proposition}
\begin{proof}
We have seen in equation \ref{eqn:Blue}
that the vanishing of the microlocal invariant
$f$ forces the vanishing of $g$, and vice-versa,
and that vanishing of both forces an elliptic
line congruence to be a Riemann sphere (in other
words, flat).
Applying our identification of the elliptic
line congruence invariants on a pseudocomplex
structure, the result is immediate.
\end{proof}

\section{Approximation by complex structures}

\begin{proposition}
\label{prop:Approx}
Suppose that $E \subset \Gro{2}{TM}$
is a pseudocomplex structure. Take $E_0 \subset \Gro{2}{TM_0}$
a complex structure (to be thought of as a flat pseudocomplex structure).
Pick points $e \in E$ and $e_0 \in E_0$, and
let $m \in M$ and $m_0 \in M_0$ be the projections
of $e, e_0$ to $M, M_0$. There is a local
diffeomorphism $f : U \to U_0$ of a neighborhood
of $m$ to a neighborhood of $m_0$ so that
$f_* e = e_0$ and so that near $e$
\[
\begin{pmatrix}
\theta_0 \\
\bar \theta_0 \\
\omega_0 \\
\bar \omega_0 \\
\pi_0 \\
\bar \pi_0
\end{pmatrix}
=
\begin{pmatrix}
1 & a_1 & 0 & a_2 & 0 & 0 \\
\bar a_1 & 1 & \bar a_2 & 0 & 0 & 0 \\
0 & b_1 & 1 & b_2 & 0 & 0 \\
\bar b_1 & 0 & \bar b_2 & 1 & 0 & 0 \\
0 & c_1 & 0 & c_2 & 1 & a_2 \bar S_{\bar 2} \\
\bar c_1 & 0 & \bar c_2 & 0 & \bar a_2 S_{\bar 2} & 1 
\end{pmatrix}
\begin{pmatrix}
\theta \\
\bar \theta \\
\omega \\
\bar \omega \\
\pi \\
\bar \pi
\end{pmatrix}
\]
where $\theta_0,$ etc. are the soldering 1-forms
of $E_0$, and so that we have
1-forms satisfying  the structure equations (i.e.\
sections of the adapted coframe bundles), which
at $e$ satisfy
\[
\begin{pmatrix}
\alpha_0 & 0 & 0 \\
\beta_0 & \gamma_0 & 0 \\
\delta_0 & \varepsilon_0 & \alpha_0 - \gamma_0 
\end{pmatrix}
=
\begin{pmatrix}
\alpha & 0 & 0 \\
\beta & \gamma & 0 \\
\delta & \varepsilon & \alpha - \gamma 
\end{pmatrix}
\]
and
\[
d 
\begin{pmatrix}
a_1 \\
a_2 \\
b_1 \\
b_2 \\
c_1 \\
c_2
\end{pmatrix}
=
\begin{pmatrix}
- \tau_1 \\
0 \\
- \tau_2 + S_{\bar 1} \pi \\
S_{\bar 2} \pi \\
- \tau_3 \\
0
\end{pmatrix}.
\]
\end{proposition}
\begin{proof} The diffeomorphism $f$ is constructed
first by taking any complex valued coordinates
$z,w$ on $M$, for which $e$ is a complex line,
and then using elementary coordinate changes
and changes of $\theta, \omega, \pi$ adapted coframing
to arrange the required equations.

To arrange the remaining equations,
which are only required to hold at $e$, we write
out an exterior differential system for the first
equations (which we have already solved), prolong it,
check that it is torsion free, and see thereby that
there is no obstruction to solving the prolongations
of all orders, at least
at any required order (but not solving it as a PDE
system---we are only solving at, say, second
order). The system will look like the above
equations together with
equations like
\[
\alpha_0 - \alpha + a_1 \bar \tau_1 + a_2 
  \left ( 
    \bar \tau_2 - \bar S_{\bar 1} \bar \pi
  \right ) = p_{11} \theta + p_{1 \bar 1} \bar \theta + p_{1 \bar 2} 
\bar \omega
\]
and so on (quite complicated), describing relations
between the various 1-forms. These coefficients
$p_{ij}$ can be chosen at will, and of course we
choose them to vanish.
Then we take a solution at whichever order
we like, arranging it to solve all the required
equations at $e$, arranging all of these arbitrary
coefficients of the prolongations to vanish at $e$,
and then repeat the manipulations
we used previously to obtain new equations 
which do not disturb the equations
at $e$, as is easy to see.
\end{proof}

\begin{corollary}
Every point
$e$ of $E$ lies in a neighborhood on which 
there are local complex coordinates $z,w$ on $M$ and 
$p$ on $E$
and adapted 1-forms $\theta, \omega, \pi$ so that
\begin{eqnarray*}
dw - p \, dz - f \, d \bar z & = & \theta + a \bar \theta \\
dz & = & b \bar \theta + \omega + c \bar \omega \\
dp & = & \pi + e_1 \theta + e_{\bar 1} \bar \theta 
+ e_2 \omega + e_{\bar 2} \bar \omega + e_3 \pi + e_{\bar 3} \bar \pi 
\end{eqnarray*}
where $a=b=c=e_{\bullet}=f=0$ at $e_0$. 

  In particular, the adapted coframes of $E$ match
those of the complex structure $E_0$ given taking the
complex coordinates $z,w,p$ to be holomorphic,
and thus $E$ and $E_0$ have the same almost complex structure
$J$ and same 4-plane field $\Theta$ at the point
$e=e_0$.
\end{corollary}

\begin{corollary} 
\label{cor:AdaptedCoords}
Take $E \subset \Gro{2}{TM}$ 
a pseudocomplex structure. Take $e \in E$ and 
let $m \in M$ be the projection of $e$ to $M$. 
There are complex valued coordinates $z,w : M \to \C{}$
defined in a neighborhood of $m$, vanishing
at $m$, and coordinates $z,w,p,q$ on $\Gro{2}{TM}$
near $e$ defined by
\[
dw = p \, dz + q \, d \bar z
\]
with $e$ represented by $z=w=p=q=0$, and 
so that $E$ is defined by the equation
\[
q = Q(z,\bar z, w, \bar w, p, \bar p)
\]
where $Q$ is a smooth function of three
complex variables so that $Q = dQ = 0$
at $e$. In particular, $E$ curves near $m$ are precisely
solutions of the equation
\[
\pd{w}{\bar z} = Q \left ( z, \bar z, w, \bar w, \pd{w}{z}, 
\overline{\left ( \pd{w}{z} \right )} \right ) \, .
\]
\end{corollary}

\section{Analysis of $E$ curves and compactness}

To allow singularities in $E$ curves, 
and study their limiting behavior,
it is natural to consider them as
parameterized rather than as submanifolds.
This is because any singular Riemann
surface can be parameterized by a (canonically
chosen) smooth Riemann surface.
We will say that a smooth map \( \phi : C \to E \)
from a Riemann surface $C$ (possibly
with boundary) to a pseudocomplex 
structure $E \subset \Gro{2}{TM}$
is a \emph{parameterized $E$ curve}
in $M$
if (1) it is $J$-holomorphic and (2) every
1-form on $E$ which vanishes on $\Theta$
pulls back through $\phi$ to $0$. 
Usually, but not always, we ask
our Riemann surface to be compact.
We will say
that our parameterized $E$ curve 
is \emph{basic} if $\phi$ is 
injective on a dense open set
and intersects each fiber of $E \to M$
on a discrete set of points. 
Henceforth, when we use the term
\emph{$E$ curve} we will mean a
parameterized $E$ curve, unless
we say otherwise.

We will also allow the selection of 
distinct marked points on the
Riemann surface $C$, and let
the Riemann surface have
ordinary double points.
(The definition of parameterized
$E$ curve from a singular Riemann
surface is just a continuous
map which lifts to a parameterized 
$E$ curve from the universal
cover.)
If $S \subset M$ is a totally real surface,
we can also consider parameterized
$E$ curves with boundary in $S$,
by which we mean that $p \circ \phi$ 
takes $\partial C$ to $S$.
A \emph{symmetry} of a
parameterized $E$ curve
is a map $C \to C$ under which 
all marked points are fixed and
$\phi$ is invariant. A parameterized $E$
curve is 
called \emph{stable} if its symmetry group is finite.

We have to be careful about parameterized
$E$ curves: they should
be thought of as curves in $M$, not in $E$,
despite the definition. However, there
is one subtlety we should keep in mind:
consider the example of $M = \CP{2}$,
with $E$ the standard complex structure.
Take the $E$ curve $C_{\varepsilon}$ defined
in affine coordinates $z,w$ on $\CP{2}$ 
and an affine coordinate $\sigma$ on $\CP{1}$ 
by
\begin{align*}
z &= \sigma + \frac{\varepsilon}{\sigma} \\
w &= \sigma - \frac{\varepsilon}{\sigma}
\end{align*}
which satisfies
\[
z^2 - w^2 = 4 \varepsilon \, .
\]
We are naturally led to suspect that
the limit as $\varepsilon \to 0$ should
be a pair of transverse lines.
The resulting curve in $E$ is
\begin{align*}
z &= \sigma + \frac{\varepsilon}{\sigma} \\
w &= \sigma - \frac{\varepsilon}{\sigma} \\
p &= \frac{\sigma^2 + \varepsilon}{\sigma^2-\varepsilon}
\end{align*}
or
\[
p = \frac{z}{w}.
\]
Naively treating $z,w,p$ as affine coordinates
on $\CP{3}$, using the equations
\[
z^2 - w^2 = 4 \varepsilon, wp=z
\]
we find that the limiting
object as $\varepsilon \to 0$ 
must be composed of four lines, not two.
This is because we traverse the exceptional
divisor $z=w=0$ twice. So $E$ curves
considered inside $E$ are slightly
different from what their
images look like in $M$.

First we will develop the local analysis.

\begin{theorem}[Elliptic regularity]
Let  \( E \subset \Gro{2}{TM} \) be a 
pseudocomplex structure.
Define a \emph{weak $E$ curve} to 
be a parameterized $E$ curve 
\[ \phi : C \to E \]
so that $\phi$
belongs
to the Sobolev space of maps
whose 1-jet is locally square integrable.
Then $\phi$ is smooth in the interior
of $C$. 
\end{theorem}
\begin{proof} This is identical to the
proof in \cite{Ye:1994}.
\end{proof}

\begin{corollary} A weak $E$ curve
has smooth prolongations to all orders,
forming pseudoholomorphic curves
in the almost complex manifolds $E^{(k)}$.
\end{corollary}

\begin{proposition} Every $E$ curve
has a prolongation to some pseudoholomorphic
curve in $E^{(k)}$ with no branch points---an immersed pseudoholomorphic curve.
\end{proposition}
\begin{proof} As in the theory of almost
complex manifolds, we have holomorphic
polynomials as lowest order terms, in
suitable complex coordinates (see 
\cite{McDuff:1992}). Moreover, the
prolongations force the order of these
polynomials down at each step, eventually
reaching order one.
\end{proof}

\begin{theorem}[Uniqueness of continuation]
Any two parameterized $E$ curves
\( \phi, \psi : C \to E \) with the same
infinite jet at a point are identical throughout
the component of $C$ containing that point.
\end{theorem}
\begin{proof} A simple application of Aronszajn's
lemma.
\end{proof}

\begin{theorem} Take a parameterized $E$ curve
\( \phi : C \to E  \) and a point $p_0 \in C$. 
Let $z,w,p$ be coordinates
as in proposition  \ref{prop:Approx}
so that $\phi(p_0)$ is the origin
of these coordinates.
Then there is a neighborhood $U$ of $p_0$ in $C$
and a holomorphic map $f=(Z,W,P) : U \to \mathbb{C}^3$
so that in the coordinates $z,w,p$ the
maps $f$ and $\phi$ have the same leading order
terms in their Taylor expansions. In particular,
the projection to $M$ in these coordinates
is holomorphic up to leading order terms.
\end{theorem}
\begin{proof} See \cite{McDuff:1992}. The basic
idea is as follows. First write out the condition on $\phi$
being $J$-holomorphic in a local coordinate
$\sigma$ on $C$ and adapted coordinates $(z,w,p)$ on $E$ 
in the form 
\[
\pd{\phi}{\bar{\sigma}}
=
M
\overline{
\left ( 
\pd{\phi}{\sigma}
\right )
}
\]
where $M$ is a complex $3 \times 3$ matrix,
vanishing at $e_0$.
This is possible because the adapted complex 
coordinates impose a complex structure
which agrees at $e_0 = \phi(p_0)$ with
the almost complex structure $J$ on $E$.
We apply Cauchy's theorem to see that
in any disk about $p_0 = 0$ in $C$:
\[
\phi
+
\frac{1}{2 \pi \sqrt{-1}}
\int_{D} 
\frac{1}{\zeta - \sigma}
M \overline{
\left ( \pd{\phi}{\sigma}
\right )}
\, d \zeta \wedge d \bar \zeta
= 
f
\]
a triple of holomorphic functions.
These intuitively represent the holomorphic
part of $\phi$ in these coordinates.
One needs to show that the map $\phi \mapsto f$
given in this way is a smooth map of
appropriate Banach spaces if the disk
$D$ is made small enough since $\phi - f$
is quite small in appropriate norms.
\end{proof}

\begin{proposition} Take $E \subset \Gro{2}{TM}$
a pseudocomplex structure and $\phi : C \to E$
a basic $E$ curve. Define the \emph{critical points} of $\phi$ to be
the points $z \in C$ where $(p \circ \phi)'(z) = 0$.
The set of critical points is a discrete subset of $C$,
except on components of $C$ that are mapped to a single point.
\end{proposition}
\begin{proof} Choose any adapted coframing $\theta, \omega, \pi$.
Now pull them back with $\phi$. By definition,
you get $\theta=0$ and $\omega$ and $\pi$ are $(1,0)$
forms, say
\[
\begin{pmatrix}
\omega \\
\pi
\end{pmatrix}
= 
\begin{pmatrix}
f \\
g 
\end{pmatrix}
d \zeta 
\]
where $\zeta$ is a local holomorphic coordinate
on $C$. Then taking differential, we find
\[
0 = \pd{f}{\bar \zeta} + f \gamma^{(0,1)} - g S_{\bar 2} \bar f
\]
where
\begin{align*}
\gamma &= \gamma^{(1,0)} d \zeta + \gamma^{(0,1)} d \bar \zeta \\
\sigma &= S_{\bar 2} \bar \omega + S_{\bar 3} \bar \pi \; .
\end{align*}
Therefore $f$ satisfies a linear first order
determined elliptic equation with smooth coefficients,
and we can apply Aronszajn's lemma to show that if
$f$ vanishes to infinite order at a point, then
it vanishes everywhere. But points where $f=0$ are
precisely points where $\phi$ is  tangent to the
fibers of $E \to M$ since 
$\theta, \bar \theta, \omega, \bar \omega$
span the semibasic 1-forms for this map and $\theta=\bar \theta =0$.
Thus $f$ can only have a discrete set
of zeros, and the set of critical points is discrete.
\end{proof}

\begin{corollary} Suppose that two parameterized basic
$E$ curves $\phi_0, \phi_1 : C \to E$ 
have projections $p \circ \phi_0, p \circ \phi_1 : C \to M$
with the same infinite jet at a noncritical point $z \in C$.
Then $\phi_0 = \phi_1$ on the connected component
of $z$ in $C$.
\end{corollary}
\begin{proof}  The equation for the
lift of the projection
\[
(p \circ \phi)^{\wedge} = \phi
\]
is easy to prove, and holds
except at critical points. Therefore,
matching of $p \circ \phi_0$ and $p \circ \phi_1$
to all orders at $z$ implies matching of $\phi_0$ and
$\phi_1$ to all orders at $z$, and therefore
we can apply uniqueness of jets of pseudoholomorphic
curves in almost complex manifolds.
\end{proof}

\begin{corollary} The same conclusion holds 
for basic parameterized $E$ curves even at a critical point. 
\end{corollary}
\begin{proof} The projections must be
asymptotic to all orders near the critical
point $z$. This implies that the
prolongations of any order of the projections
must be asymptotic to all orders as we approach
the critical point. As above
\[
\left  ( p \circ \phi_j \right )^{\wedge} = \phi_j
\]
at noncritical points. This
implies, because the critical points are discrete,
that the $\phi_j$ must be asymptotic at $z$ to all
orders. Therefore, they must agree by
Aronszajn's lemma.
\end{proof}

\begin{lemma} Suppose that $C \subset M$
is an embedded smooth $E$ curve, and $m \in C$
is a point of $C$. There are local coordinates
$z,w,p$ on $E$ near $m$, with $z,w$ defined
on $M$, so that $C$ is cut out near $m$ by
the equations
\[
w = p = 0 \; .
\]
Moreover the is an adapted coframing \( \theta, \omega, \pi \)
defined near $T_m C$ in $E$ so that
\[
\begin{pmatrix}
\theta \\
\omega \\
\pi 
\end{pmatrix}
= 
\begin{pmatrix}
dw \\
dz \\
dp
\end{pmatrix}
\]
at all points of $C$ near $m$.
\end{lemma}
\begin{proof} The relevant equivalence
problem concerns the pullback of $B \to E$
to the lift of $C$, $\tilde{C} \subset E$.
We can further adapted the coframes
$\theta_0, \omega_0, \pi_0$ 
of this pullback by asking that 
$\pi_0, \bar \pi_0$ vanish on tangent spaces of $C$.
This forces structure equations
\[
\theta = \bar \theta = \pi = \bar \pi = 0
\] 
which gives us as structure equations
on this bundle
\[
d \omega = - \gamma \wedge \omega
\]
which are the structure equations
of a complex structure on a surface,
and so local equivalence with the 
flat example follows by the Newlander--Nirenberg
theorem. Then $z,w,p$ can be arbitrarily extended
off of $\tilde{C}$ so that $z,w$ are functions
on $E$.
\end{proof}

\begin{proposition} Let $E \subset \Gro{2}{TM}$ 
be a proper pseudocomplex structure.
Suppose that $\phi_0 : C_0 \to E$
and $\phi_1 : C_1 \to E$ are $E$ curves, 
and that there are convergent sequences of points
\[
s_j \to s \in C_0, \quad t_j  \to t \in C_1
\]
so that
\[
p \circ \phi_0 \left ( s_j \right ) = p \circ \phi_1 
\left ( t_j \right ) \; .
\]
Also suppose that $\phi_0'(s) \ne 0$.
Then there is a holomorphic map $f$ of 
a neighborhood of $s$ to a neighborhood of $t$
so that
\[
\phi_1 = \phi_0 \circ f \; .
\]
\end{proposition}
\begin{proof}
If we were to try to prove this result 
by imitating the almost complex case, we would
simply look in coordinates
of the type guaranteed by the previous lemma,
so that $C_0$ is cut out by $w=p=0$. We can
use $z$ as a local holomorphic coordinate
on $C_0$, and arrange that $z=0$
is the point $s$. Take $\zeta$ 
a local holomorphic coordinate on $C_1$,
so that $\zeta=0$ is the point $t$. Now
we have $C_1$ described by functions
\[
z(\zeta), w(\zeta), p(\zeta)
\]
with $w(\zeta)$ having infinitely many
zeros near $\zeta=0$. The lowest order
terms of 
\[
\left ( z(\zeta), w(\zeta), p(\zeta) \right )
\]
at $\zeta=0$
must be holomorphic polynomials, unless
they vanish to all orders. In the almost
complex case, this actually determines
that $w(\zeta)$ itself has holomorphic
lowest order term. But then,
$w(\zeta)$ cannot vanish at infinitely many
points approaching $\zeta=0$, unless
this lowest order term vanishes,
because it will dominate:
\[
w(\zeta) = w_0 \zeta^k + \ldots
\]
would give $w(\zeta) \ne 0$ near
$\zeta = 0$. Therefore $w(\zeta)$
vanishes to infinite order at $\zeta = 0$.
Because of the equation
\[
dw = p \, dz + q(z,w,p) \, d \bar z
\]
we have $p(\zeta)$ and $q(z(\zeta),w(\zeta),p(\zeta))$
vanishing to all orders in $\zeta$.
Now this gives us
\[
\omega = dz, \pi = dp
\]
at $\zeta = 0$, up to infinite order,
and since $\omega$ is a $(1,0)$ form, this
tells us that $z(\zeta)$ is holomorphic
to all orders. We can  
$\omega \wedge d \zeta = 0$
to show that the formal holomorphic series given by
the Taylor expansion of $z(\zeta)$ must converge
(because the terms are controlled, by differentiating
this equation, in terms of values of $z,w,p$).
The result then follows by uniqueness of continuation,
using Aronszajn's lemma.

  However, in the pseudocomplex case, it
is unclear that $w(\zeta)$ must have holomorphic
lowest order term, simply because $(z(\zeta),w(\zeta),p(\zeta))$
does; it might look like
\[
p(\zeta) = p_0 + \ldots
\]
with $p_0 \ne 0$, and then we would have lowest
order terms (constant terms) of $(z,w,p)(\zeta)$
holomorphic. The problem is essentially
to show that $p(0)=0$, so that in our coordinates
the two curves strike in $E$.

  But this is not too hard: each real line
at a point of $m$ is contained in a unique $E$ plane
in $E_m \subset T_m M$. Now take each pair of 
points $\phi_0 (s_j), \phi_0 (s_k)$ and draw a line
between them in some local coordinates. Then
the $E$ planes containing these lines must
converge to $(p \circ \phi_0) ' (s) \cdot T_s C_0$,
the $E$ plane tangent to $C_0$ at $s$.
But then, we could have used points close to
the $s_j$ instead of the actual $s_j$, or
points close to the $t_j$ as well. We can
pick points of $C_1$ which are not critical 
for $\phi_1$, and carry out the same construction.
This shows that the values of $p(\zeta)$ must tend
to $0$ as $\zeta \to 0$. Therefore $C_0$ and $C_1$
have the same osculating complex geometry at
$m$, so their lowest order terms are holomorphic
there. The story is now the same as in the
almost complex case.
\end{proof}

\begin{theorem} Two $E$ curves have a discrete set
of intersection points in $M$ (finite, if they are compact).
Each intersection point, after small perturbation,
becomes a positive number of transverse intersections.
Moreover, the local intersection number at
an intersection is one precisely when the intersection
is transverse.
\end{theorem}
\begin{proof} The proof is essentially as in
\cite{McDuff:1994}. 
\end{proof}

\begin{definition} Take $M$ a 4-manifold with
Riemannian metric $g$ and \( E \subset \Gro{2}{TM} \)
a pseudocomplex structure for $M$. Let
$\Omega \in \nForms{2}{M}$ be a symplectic
form on $M$. The \emph{weight} of $E$
with respect to $g$ and $\Omega$ is
the smallest constant $a$ so that
\[
a \cdot \Omega(u,v) \ge 1
\]
for $u,v \in T_m M$ any oriented $g$ orthonormal basis
of a 2-plane $\textit{span}\{ u,v \}$ belonging to $E$.
The \emph{weight} of $\Omega$ with respect to $g$ is
the smallest constant $b$ so that
\[
\Omega(u,v) \le b
\]
whenever $u,v \in T_m M$ are orthogonal
$g$-unit vectors. The \emph{width} of $E$ with
respect to $g$ and $\Omega$ is the
smallest  constant $c$ so that all of the osculating
complex structures $J_e : T_m M \to T_m M$
for $e \in E$ satisfy
\[
\left | J_e v \right | \le c
\]
for $v \in T_m M$ any $g$-unit vector.
\end{definition}

\begin{proposition} Suppose that $M$ is a 
4-manifold, that $\Omega \in \nForms{2}{M}$
is a symplectic form, and that \( E \subset \Gro{2}{TM} \)
is a pseudocomplex structure. Let $g$ be 
a flat Riemannian metric. Let
$a$ be the weight of $E$, $b$ the
weight of $\Omega$ and $c$ the width
of $E$. Suppose that $a, b, c > 0$ and that
$\Omega = d \vartheta$
is exact, for some 1-form $\vartheta$.
Take any parameterized $E$ curve 
\( \phi : C \to E \) and any Borel 
set \( B \subset M \) and rectifiable
current $X$ on $M$ with 
\[
\partial X = \partial \left ( p \circ \phi(C) \rhook B \right )
\]
\[
\textit{Mass} \left ( p \circ \phi(C) \rfloor B \right )
\le abc \, \textit{Mass} \left ( X \rfloor B \right ).
\] 
\end{proposition}
\begin{proof} See \cite{Bangert:1998}.
\end{proof}

\begin{proposition} Suppose that
$E, \Omega, g$ and $\phi : C \to E$ 
are as in the
last proposition, with $a, b, c > 0$. Suppose 
further that $g$
is a flat Euclidean metric, and that $C$ is 
compact with boundary, and that
\(
p \circ \phi (\partial C)
\)
has support (as a current) contained in a 
ball of radius $r$. Then for any integer $j > 0$
and real number $\varepsilon > 0$,
the support (as a current) of $p \circ \phi(C)$ is
contained in a ball of radius
\[
\varepsilon^j r +
\left ( \frac{
	1 - \frac{1}{abc}}{1-\varepsilon^{-2}} \right )^{j/2}
\sqrt{
\frac{abc}{4 \pi}
\textit{Mass} \left (
S \rfloor \left ( \R{2n} \backslash B(r) \right )
\right )
}.
\]
\end{proposition}
\begin{proof} See \cite{Bangert:1998}.
\end{proof}

\begin{theorem}[Removable singularities]
Let $E \subset \Gro{2}{TM}$ be a proper pseudocomplex
structure. Suppose that $\phi : C \to E$ is an $E$ curve,
with $C$ a punctured Riemann surface, and the
image of $p \circ \phi$ is contained in a compact subset 
of $M$. Then $\phi$
extends across the puncture precisely when 
the area of the image of $\phi$ is finite.
\end{theorem}
\begin{proof} As in \cite{Gromov:1985}.
\end{proof}

\begin{theorem}[Gromov compactness]
Take \( E_j \subset \Gro{2}{TM} \) a 
sequence of proper
pseudocomplex structures, $S_j \subset M$ a
sequence of
surfaces, with $S_j$ totally real for $E_j$ 
and a sequence
of stable parameterized $E$ curves 
\( \phi_j : \left ( C_j , \partial C_j \right ) \to 
\left ( E, p^{-1} S_j \right ) . \)
Suppose that $E_j \to E$ converges
to a proper pseudocomplex structure,
and $S_j \to S \subset M$ converges to a totally
real surface for $E$. Suppose that there
are Riemannian metrics on the $E_j$ 
converging to one on $E$ for which
the area of the image of any of the $\phi_j$
is less than some bound $A$, and that
the $p \circ \phi_j$ have images contained inside
some compact set $K \subset M$.
Then there is some subsequence $\phi_{j_k}$ 
and a unique stable parameterized $E$ curve $\phi : C \to E$
so that \( \phi_{j_k} \) converges
to $\phi$ in the Gromov--Hausdorff sense.
\end{theorem}
\begin{proof} This follows from any of the
myriad proofs of Gromov compactness, by
applying the canonical almost complex structure
on the $E_j$.
\end{proof}

\begin{corollary}
\label{cor:UnifCon}
Uniform convergence on compact sets for
a sequence of parameterized $E$ curves
implies uniform convergence on compact
sets for all derivatives.
\end{corollary}
\begin{proof} Use prolongation to all orders.
\end{proof}

\begin{proposition} 
\label{prop:CmptSym}
Uniform convergence on compact sets of
a sequence of symmetries of a pseudocomplex
4-manifold implies uniform convergence on
compact sets of all derivatives.
\end{proposition}
\begin{proof} Apply the symmetry to
various $E$ curves.
\end{proof}

Consider a pseudocomplex structure in
local adapted coordinates $(z,w,p)$. It
looks like a function
\[
q(z,w,p)
\]
which is smooth, and satisfies
\[
0 = q(0) = dq(0) \; .
\]
A parameterized $E$ curve with local
holomorphic parameter $\sigma$ must satisfy
\[
dw = p \, dz + q(z,w,p) \, d \bar z
\]
and also have 
\[
\omega \wedge d \sigma = \pi \wedge d \sigma = 0 \; .
\]
This can be written in our local coordinates
as equations
\begin{equation}
\label{eqn:Stuffy}
\begin{pmatrix}
\pd{z}{\bar \sigma} \\
\pd{p}{\bar \sigma} 
\end{pmatrix}
=
\begin{pmatrix}
a_1 & a_2 & a_3 \\
b_1 & b_2 & b_3 
\end{pmatrix}
\begin{pmatrix}
\pd{z}{\sigma} \\
\pd{w}{\sigma} \\
\pd{p}{\sigma}
\end{pmatrix}^{\dag}
\end{equation}
where the ${}^{\dag}$ indicates complex conjugation,
and $a_j, b_k$ depend on $z,w,p$ smoothly,
algebraically determined from $q(z,w,p)$ and its 
first derivatives. We can try to produce holomorphic
approximations to $z,p$ by writing
\[
\begin{pmatrix}
Z \\
P 
\end{pmatrix}
(\sigma)
=
S_U \left [ \begin{matrix} z \\ p \end{matrix} \right ] (\sigma)
=
\frac{1}{2 \pi \sqrt{-1}}
\int_{\partial U}
\begin{pmatrix}
z(\zeta) \\
p(\zeta)
\end{pmatrix}
\frac{d \zeta}{\zeta - \sigma}
\]
for some small region $U$ around the origin,
so that
\[
\begin{pmatrix}
z \\
p
\end{pmatrix}
=
\begin{pmatrix}
Z \\
P
\end{pmatrix}
+
T_U \begin{pmatrix}
\bar \partial z \\
\bar \partial p
\end{pmatrix}
\]
where
\[
T_U \left [ g(\zeta) d \bar \zeta \right ](\sigma)
=
\frac{1}{2 \pi \sqrt{-1}} \int_U
\frac{g(\zeta) d \bar \zeta \wedge d \zeta}{\sigma - \zeta} \; .
\]

Now suppose that we start with some
arbitrary complex valued functions $z,p$ of
$\sigma$, and determine a function $w$
by solving the ODE system
\[
dw = p \cdot dz + q(z,w,p) \, d \bar z
\]
radially away from the origin,
in a small open set in the $\sigma$ complex
plane,
and then construct $Z,P$ by
\[
\begin{pmatrix}
Z \\
P
\end{pmatrix}
=
\begin{pmatrix}
z \\
p
\end{pmatrix}
-
T_U 
\begin{pmatrix}
\bar \partial z \\
\bar \partial p
\end{pmatrix}
\]
but we plug in the expressions
from the right-hand side of equation
\ref{eqn:Stuffy} instead of the partials
\( \bar \partial z, \bar \partial p. \)
This gives a map from arbitrary functions
$z,p$ defined near 0, satisfying
some appropriate smoothness condition,
to functions $Z,P$ satisfying a similar
condition. Using the same arguments
as in \cite{McDuff:1992}, we find
that this map is invertible near 0,
on a suitable function space, for small
enough open set $U$,
and that the holomorphic choices of $Z(\sigma),P(\sigma)$
correspond precisely to the $z(\sigma),p(\sigma)$
so that the triple $z(\sigma),w(\sigma),p(\sigma)$ is 
an $E$ curve.

\begin{proposition} Let $E_t \subset \Gro{2}{TM}$
be a family of pseudocomplex structures, parameterized
by a real variable $t$. For any $E_0$ curve $\phi : C \to E_0$,
there is a neighborhood $U \subset C$ of any point
and a
deformation $\phi_t : U \to E_t$
of $\phi$ which is smooth in $t$, defined for some
open set of $t$ values near $0$, giving an $E_t$
curve for each $t$.
\end{proposition}
\begin{proof} The proof is due to McDuff \cite{McDuff:1992}.
The idea is to carry out the above construction
of $Z,P$ for the initial curve, which is governed
by the function $q_0(z,w,p)$, say. Then we invert
the construction for the functions $q_t(z,w,p)$
corresponding to each $E_t$. This can only be done
locally.
\end{proof}

\begin{proposition} Given any $E$ curve
$\phi : C \to E$ we can make a small deformation
of $E$, say $E_t$, and of $\phi$, say $\phi_t$, 
an $E_t$ curve, so that $E_0 = E$ and $E_1$ is
flat (complex) near all singular values of $\phi_1$.
\end{proposition}

\section{Dual curves}

\begin{theorem} Let $E$ be a pseudocomplex structure on
$\CP{2}$, tamed by a symplectic structure. Then
the moduli space of $E$ spheres representing
the hyperplane homology class is smooth and
compact and diffeomorphic to $\CP{2}$. 
\end{theorem}
\begin{proof} By results of Taubes, there is
only one symplectic structure on $\CP{2}$ up
to rescaling, so that we can suppose that
our symplectic structure is the usual one on
$\CP{2}$, and from here the rest of the
arguments are the same as in Gromov's \cite{Gromov:1985}.
\end{proof}

  Now suppose that $E$ is a pseudocomplex structure
on a 4-manifold $X$, and that $Y$ is a 4 parameter
family of $E$ curves. Let $Z$ be the set of 
pointed $E$ curves whose underlying $E$ curves
belong to the family $Y$. We have maps
\[
\xymatrix{
& Z \ar[dl]_{\lambda} \ar[dr]^{\rho} & \\
X & & Y
}
\]
with $\lambda$ given by forgetting the $E$ curve,
and $\rho$ by forgetting the point. We also
have a map
\[
(p,C) \in Z \to T_p C \in E \; .
\]
Suppose that this map is a local diffeomorphism
(as it is in the case of $\CP{2}$, essentially
proved in \cite{Gromov:1985}). Then we can
pull our $G$ structure from $E$ back to $Z$,
and reduce the structure group by taking only
the adapted coframings $\theta_0, \omega_0, \pi_0$
for which the fibers of $Z \to Y$ satisfy
\[
\theta = \bar \theta = \pi = \bar \pi = 0 \; .
\]
This reduces the structure group to the group
of complex matrices
\[
\begin{pmatrix}
a & 0 & 0 \\
b & c & 0 \\
d & 0 & ac^{-1}
\end{pmatrix}
\]
so that $a,c \ne 0$. This reduction forces the
1-form $\epsilon$ to be semibasic, and we 
can arrange
\[
\epsilon = E_1 \bar \theta + E_2 \bar \pi \; .
\]
Now define 1-forms on this principal bundle by
\[
\begin{pmatrix}
\theta' \\
\omega' \\
\pi' 
\end{pmatrix}
=
\begin{pmatrix}
\theta \\
\pi \\
- \omega
\end{pmatrix}
\qquad
\begin{pmatrix}
\alpha' & 0 & 0 \\
\beta' & \gamma' & 0 \\
\delta' & \epsilon' & \alpha' - \gamma'
\end{pmatrix}
=
\begin{pmatrix}
\alpha & 0 & 0 \\
\delta & \alpha-\gamma & 0 \\
- \beta & \sigma & \gamma
\end{pmatrix}
\]
and define complex valued functions according to
\[
\begin{pmatrix}
S_1' \\
S_2' 
\end{pmatrix}
=
\begin{pmatrix}
E_1 \\
E_2
\end{pmatrix}
\qquad
\begin{pmatrix}
T_2' & T_3' \\
U_2' & U_3' \\
V_2' & V_3'
\end{pmatrix}
=
\begin{pmatrix}
-T_3 & T_2 \\
V_3 & - V_2 \\
-U_3 & U_2
\end{pmatrix}
\]
Then these new 1-forms and functions satisfy
the structure equations \ref{eqn:struc}.
Therefore, by Proposition \ref{prp:Rebuild},
this is a pseudocomplex structure, at least
locally, and the map $Z \to Y$ presents
$Y$ as the base space of this structure,
so that this is locally a pseudocomplex
structure on $Y$ (it could be ``multivalued'').
Moreover, the $E$ curves for $Z \to Y$
are the same (as surfaces in $Z$) as the
$E$ curves of $Z \to X$.

\begin{theorem} Suppose that the
pseudocomplex structure given (locally) by $Z$
is an almost complex structure (locally) both
on $X$ and on $Y$. Then it is a complex structure
on both.
\end{theorem}
\begin{proof} This follows immediately from
the identification of the equations on 
microlocal invariants of almost complex structures:
we would need 
\[
S_1 = S_2 = T_3 = U_3 = 
S_1' = S_2' = T_3' = U_3' = 0
\]
which (when we unwind these equations
and differentiate the structure equations)
forces all of our invariants to vanish; hence local
equivalence with the flat case.
\end{proof}

\section{Gromov--McDuff--Ye intersection theory}
\label{sec:GMY}

Intersections behave as they do in complex
geometry, because of Proposition \ref{prop:Approx}.

\begin{proposition} Let \( \phi : C \to E \)
be a continuous map which is an $E$ curve
in a neighborhood of a point $p_0 \in C$.
Define the local self-intersection number
of $C$ at $p_0$ to be the local 
self-intersection number of the map
\( p \circ \phi . \) This number is
always positive unless  
$(p \circ \phi)' \left (p_0 \right )$
is an injection. The local intersections
in $M$ of two $E$ curves 
are always positive. The local intersection
number is 1 only if the curves intersect
transversely.
\end{proposition}

\subsection{Blowup and adjunction}

Take $E \subset \Gro{2}{TM}$ a proper pseudocomplex
structure. To define blowup of a point, say $m \in M$, 
first we need the analogues of lines through
that point. Any family of
$E$ curves will do, as long as it contains
one line through each point in each ``direction'',
i.e.\ each $E$ tangent line. First, by
cutting out a neighborhood of $m \in M$, 
and pasting in $\CP{2}$, following arguments
of Bangert in \cite{Bangert:1998}, we can
arrange that the resulting pseudocomplex
structure be tamed by the usual symplectic 
structure on \CP{2}
and that it agree  with the usual complex
structure outside some small neighborhood,
and we can even get it to be as close
as we like to the usual complex structure
on $\CP{2}$ using our arguments for approximation
by complex structure. Now we can employ
Gromov's compactness theorem on $E$ to
ensure that the small deformation of the
usual complex structure to $E$ preserves
the family of lines. Using the intersection
theory arguments of McDuff, given
in Section \ref{sec:GMY}, we find that
these curves must remain smooth under
such deformation, and that they intersection
in $M$ in only one place, transversely.
The problem is to get one of them
through each point in each direction.
This is not difficult: we find that
any pair of distinct points,  are joined
by a unique $E$ curve in our family, 
by arguments of Gromov \cite{Gromov:1985}.
This $E$ curve varies smoothly with the choice
of points.

  Now take a point of one of these
$E$ curves, and look at the linearization
of the $E$ operator on that curve.
One finds (for example, by index theory
of \CP{1}) that its kernel has real dimension
4, which matches the index, 
and that the space of solutions
of the linearization vanishing at a point
is precisely two-dimensional, again matching
the index for that problem. So at
least locally, the space of $E$ curves
from our family, passing through a 
given point, is a smooth surface, and 
no two such curves sufficiently near
one another can have the same tangent
plane at that point. Therefore, the
family of curves is locally identified
with their tangent planes. Properness
ensures that this is true globally
as well.

  The end result then is that locally,
we can construct such a family of
$E$ curves, looking ``like lines''.
We can then define blowup imitating
the usual definition: away from $m \in M$
define the blowup of $M$ at $m$ to be just identified
with $M$ itself; but in a small neighborhood
of $m$, we define it to consist of pairs
$(n,C)$ where $n \in M$ is a point
belonging to an $E$ curve $C$ which
belongs to our family. The exceptional
divisor is then just $E_m$.

 Unfortunately, there is in general
no pseudocomplex structure on the blowup
which (1) agrees with the usual structure
away from the exceptional divisor
(where we can identify the blowup with
$E$) and which (2) will extend continuously
across the exceptional divisor.
However, we can define the lift of
a curve $C_0$ through $m$ to the
blowup: to each point $n \in C_0$
we associate the point $n$ in the
blowup if $n$ is away from $m$,
and for $n$ nearby $m$ associate the point
of the form $(n,C)$. This will extend
to a smooth real surface in the
blowup, which we see easily as follows:
the blowup near the exceptional divisor
is embedded into $E$ by sending a point
$(n,C)$ of the blowup to $T_n C \in E$.
This map
\[
(n,C) \mapsto T_n C \in E
\]
is easily seen to be an embedding
of a neighborhood of the exceptional
divisor (for example, in adapted coordinates).
Now we can see easily that the
blown up curve inside the blowup of $M$ at $m$
sits in $E$ as a surface
asymptotic to the lift $\hat C_0$ of $C$
(the set of tangent lines to $C$). This
allows us to extend the blowup curve
across the exceptional divisor smoothly.
Intersections with the exceptional divisor
correspond (including multiplicities)
with intersections of the lift $\hat C_0 \subset E$
and the exceptional divisor $E_m$, which
we can see again by deforming to the
complex case.

  Since the blowup does not have a
pseudocomplex structure (the obvious choice
pulled back from $M$ does not extend differentiably
across the exceptional divisor, in general)
it is probably not natural to approach the
study of $E$ curves via blowup. In particular,
repeated blowup is not possible, while repeated
prolongation is. However, results on prolongations
always require differential geometry and
are not topological.

The notion of blowdown is particularly
intriguing here, but even in the
presence of an exceptional divisor,
blowdowns will not generally be well
defined.

 Totally real surfaces can be treated in two
ways: first, the way I approached them in my
thesis, where I simply took as the lift of
a surface all of the 2-planes in $E$ that
live at that point of $M$. But to treat totally
real surfaces, one can look at each 2-plane
belonging to $E$ above a point of my surface $S$,
look at each $e \in E$ above a point $m \in S$.
Then ask if $T_m S$ is totally real for 
$J_e : T_m M \to T_m M$,
i.e.\ $T_m S$ is not $J_e$ invariant, where
$J_e$ is the osculating complex structure.
Then I ask further if $e$ is invariant 
(up to reorientation) under
the complex conjugation across $T_m S$ in
$T_m M$. Recall that every totally real 2-plane
in a complex two-dimensional vector space has
a conjugation operator. These $e$ are precisely
the complexifications of real lines (they
strike $T_m S$ in a real line). So they can be 
identified with a circle bundle (the real
projectivization of $TS$) above $S$ in
the case of an almost complex manifold.
Call the total space $S''$.
But for a pseudocomplex manifold, there
is more to it since the osculating
Riemann sphere may not be constant.
$S''$ is a circle bundle over $S$,
in the proper case. Therefore, a
3-manifold. This allows much tighter
adaptation of frames. If $C$ is an $E$
curve with boundary in $S$, then its
lift is a pseudoholomorphic curve in $E$
with boundary in $S''$. 

By carrying out infinite prolongation,
we will find that any two $E$
curves with boundary, whose boundaries
have the same infinite jet at a point, 
must be identical. The
reasoning: we will find that it is possible
to complexify a real curve to each order,
by analogy with the case of complex
surfaces. No finite order obstruction
emerges since we easily see involutivity.
This determines the infinite
jet of a pseudoholomorphic curve.
Now we apply Aronszajn's lemma.

The possibility of slightly extending
smooth $E$ curves
with boundary is immediate from the
continuity method.

If $e \in S'',$ then we have $T_m S$ totally
real for $J_e : V \to V,$ so that $T_m S$
is also totally real for any approximating
complex structure which is tangent to
our pseudocomplex structure. We can therefore
use holomorphic coordinates for the
complex structure in which $S$ 
matches $\mathbb{R}^2$ to any desired
order at a single point.

The three manifold $S''$ is totally real inside $E$.
This $S''$ is the first step in an infinite
prolongation of $S$.

\section{The canonical bundle}

Define the canonical bundle $K$ of a pseudocomplex
structure $E$ on a 4-manifold $M$ to be
the complex line bundle associated to the
representation
\[
\rho 
\begin{pmatrix}
a & 0 & 0 \\
b & c & 0 \\
d & e & ac^{-1} 
\end{pmatrix}
=
ac
\]
of the structure group $G$.
Note that this line bundle lives on $E$ and
not on $M$.
For a proper pseudocomplex
structure, it follows from local deformation
to a complex structure that $K$ restricts
to each fiber of $E \to M$ to be 
a trivial holomorphic line bundle.

We can define a $\bar \partial$ operator
on $K$ by the following apparatus:
a section of $K$ is identified with a 
function $f : B \to \C{}$ satisfying
equivariance
\[
r_g^*f = \frac{f}{\rho(g)} .
\]
We define
\( \bar \nabla f \)
to be the $(0,1)$ part of the 1-form
\[
df - f \left ( \alpha + \gamma \right ) .
\]
However, there is no reason to believe
that local holomorphic sections of $K$
exist. The reader might try 
to determine when such sections exist,
as an
exercise in using the structure equations.

Now we have a canonical bundle for any
pseudocomplex structure and a canonical
bundle on any complex curve, so we can
take a parameterized $E$ curve
\[
\phi : C \to E
\]
and pull back the canonical bundle $K_M$ from $E$,
\[
\phi^*K_M
\]
and look for an adjunction formula,
comparing
\[
\phi^*K_M \text{ and } K_C \; .
\]

\begin{theorem}[Adjunction]
Let $E \subset \Gro{2}{TM}$ be a proper pseudocomplex
structure.
Let $C$ be a compact smooth Riemann surface,
and \( \phi : C \to E \) a parameterized $E$ curve.
Then
\[
- \chi_C = K_M \cdot C + C \cdot C - 2 \delta
\]
where $\delta$, the \emph{embedding defect}
is a nonnegative integer, vanishing only
if $p \circ \phi : C \to M$ is an embedding.
\end{theorem}
\begin{proof} Suppose that $p \circ \phi : C \to M$
is an embedding. Then we can see from the
structure equations that the normal
bundle to $C$ in $M$, call it \( \nu_C, \) 
is the pullback by $\phi$ of 
the vector bundle associated to the
principal bundle $B \to E$ and the representation
\[
\begin{pmatrix}
a & 0 & 0 \\
b & c & 0 \\
d & e & ac^{-1}
\end{pmatrix}
\mapsto
a^{-1}
\]
while the canonical bundle of $M$ with respect
to $E$ is associated to the representation
\[
\begin{pmatrix}
a & 0 & 0 \\
b & c & 0 \\
d & e & ac^{-1}
\end{pmatrix}
\mapsto
ac
\]
and the canonical bundle of $C$ is
associated to
\[
\begin{pmatrix}
a & 0 & 0 \\
b & c & 0 \\
d & e & ac^{-1}
\end{pmatrix}
\mapsto
c.
\]
Therefore 
\[
K_M \otimes C = K_C
\]
which gives the result immediately:
\begin{alignat*}{2}
- \chi_C &= \int_C c_1 \left ( K_C \right ) 
&& \text{by Riemann--Roch} \\
         &= \int_C c_1 \left ( K_M \otimes \nu_C \right ) 
&& \\ 
         &= \int_C c_1 \left ( K_M \right ) 
	+ \int_C c_1 \left ( \nu_C \right ) 
&& \\
         &= K_M \cdot C + C \cdot C. 
\end{alignat*}

The problem in the nonembedded case
comes from $C \cdot C$, the self
intersection number,  not being identical
to $\int_C c_1 \left ( \nu_C \right )$.
We might imitate the complex proof and carry out
blowup, but the blowup is not a pseudocomplex
manifold. To deal with this, we can follow 
McDuff's arguments, making a deformation
of pseudocomplex structure, say as $E_t$, near each singular
value of $p \circ \phi$, and deforming $\phi$
as we do, say as $\phi_t$,  so that $E_0=E$
and $E_1$ is complex near each singular
point of $\phi_1$. Then we can employ
blowup at each singular point.
\end{proof}

\section{Darboux's method of integration}

\begin{definition} A pseudocomplex structure 
$E \subset \Gro{2}{TM}$ is called \emph{Darboux
integrable} if near any point of $E$ there are
two independent holomorphic functions.
\end{definition}

We can use such functions as follows:
take $f,g : E \to \mathbb{C}$ to be 
such functions. Then any holomorphic
curve in $\mathbb{C}^2$ has preimage
under $(f,g)$ some almost complex 4-manifold
in $E$. The equations $\theta=\bar \theta=0$
restrict to this 4-manifold to form
a holonomic plane field, i.e.\ a
system of ordinary differential equations,
to find all $E$ curves which map under 
$(f,g)$ to the chosen curve.
Conversely, if we have any $E$ curve,
then $(f,g)$ must restrict to it to
form a pair of holomorphic functions
on the $E$ curve so that for generic
$E$ curves this will trace out a 
holomorphic curve in $\mathbb{C}^2$. 
Consequently, $E$ curves can be found
by solving ordinary differential equations
only.

Indeed, we need not require $f,g$ to be
holomorphic, but only to be holomorphic
on all $E$ curves; however this turns
out to force $f$ and $g$ to be holomorphic
on $E$. Indeed, plugging in to the
structure equations, using the equation
\[
d^2 f = 0
\]
and writing out
\[
df = f_1 \theta + f_2 \omega + f_3 \pi + f_{\bar 1} \bar \theta +
f_{\bar 2} \bar \omega + f_{\bar 3} \bar \pi,
\]
we find that if such
$f$ and $g$ exist, then we have linear equations
among the rows of
\[
\begin{pmatrix}
T_2 & T_3 \\
U_2 & U_3 \\
V_2 & V_3
\end{pmatrix}.
\]
All of these invariants vanish precisely when $E$ is 
a complex structure on a complex surface $M$. 
Darboux integrability requires two
functions $f$ and $g$ which have independent
complex linear differentials modulo $\bar \theta$.
So this implies that, in some adapted coframing
\[
U_2 = U_3 = V_2 = V_3 = 0 \; .
\]
There is no complete classification of the
possibilities, but if we assume that either
(1) $T_3 \ne 0$ everywhere, or (2) that both 
$T_3=0$ everywhere and $T_2 \ne 0$ everywhere,
or (3) that $T_2=T_3=0$ everywhere, then we can
split up our study into each of these cases.
Similarly, in case (2), we have to split up
into possibilities depending on whether some
other invariant vanishes; the reader will find details 
and computations in 
\cite{McKay:1999}. 
The resulting pseudocomplex structures
can be tabulated as follows (at least
locally):
\begin{enumerate}
\item
  Complex surfaces. Locally, up to diffeomorphism
  there is a unique example. The symmetry pseudogroup
  is the group of local biholomorphisms, depending
  on 4 real functions of 2 real variables.
\item \label{itm:DrbInt}
  Take $C$ a complex curve, and $\mathcal{K}$ its
  canonical bundle. Then the equation
  \[
    d \xi = \bar{\xi} \wedge \xi
  \]  
  defines an almost complex structure on the total
  space of the canonical bundle. If we take
  a local holomorphic coordinate $z$ on $C$,
  this differential equation can be written
  \[
    \pd{w}{\bar{z}} = |w|^2
  \]
  in local holomorphic coordinates.
  Locally, up to diffeomorphism, there
  is a unique example. The symmetry pseudogroup
  is the group of local biholomorphisms of
  the complex curve $C$, depending on
  2 real functions of 1 real variable.
\item \label{itm:ACC}
  An almost complex example, which 
  has coordinates $Z,W$ and can be described
  as the almost complex structure for which the
  1-forms
  \begin{align*}
     & dW -  \frac{W \bar{Z}}{1-|Z|^2} \,  dZ 
          - \frac{\bar{W}}{1-|Z|^2} \, d \bar{Z} \\
     \intertext{and}
     & dZ 
  \end{align*}
  are holomorphic. To find the $E$ curves,
  take $Z(P)$ any holomorphic function of one
  complex variable, and integrate the ordinary
  differential equation
  \[
    0 = dW - \left ( \frac{W \bar{Z}}{1-|Z|^2} + P(Z) \right ) \, dZ 
        - \frac{\bar{W}}{1-|Z|^2} d \bar{Z}.
  \]
  The manifold $M$ is a bundle of complex
  curves over a complex
  curve. The base has parameter $Z$ and a hyperbolic
  metric. The fibers 
  $Z=\text{constant}$ are complex curves,
  so this is an almost complex structure
  canonically defined on some line bundle
  over a hyperbolic complex curve. Its geometric
  meaning is unclear to the author. 
  Locally, up to diffeomorphism, there
  is a unique example.
  These $Z,W$ coordinates are determined up
  to the following transformations:
  \begin{itemize}
     \item  
        Any orientation preserving 
        hyperbolic metric isometry
        of the curve parameterized by the 
        $Z$ variable. Under these
        transformations, the $W$ variable
        behaves as a spinor.
     \item
        Take $f(Z)$ any holomorphic function,
        and write down the complex valued
        functions \( a \left ( Z,\bar{Z} \right ), 
           r \left ( Z,\bar{Z} \right ) \) 
        determined by the equations
        \begin{align*}
           r \left ( Z,\bar{Z} \right ) &= 
                Z^2 f(Z) +
                 \overline
                 {
                   \frac{d}{dZ} Z^2 f(Z) 
                 } \\
           a \left ( Z, \bar{Z} \right ) &= 
                 \frac{\bar{r} + r \bar{Z}}{1-|Z|^2}.
        \end{align*} 
        Then define new coordinates $\hat{Z}, \hat{W}$ by
        \begin{align*}
           \hat{Z} &= Z \\
           \hat{W} &= W + a \left ( Z, \bar{Z} \right ).
        \end{align*}
     \end{itemize}
\item
  Take $F \left ( z, \bar{z}, w, \bar{w}, p, \bar{p} \right )$
  any solution of the elliptic determined quasilinear
  system of partial differential equations
  \[
    0 = \pd{F}{\bar{p}} + \bar{F} \pd{F}{\bar{z}}
        + \frac{\bar{w}}{1-|p|^2} \pd{F}{w} 
        + \frac{p \bar{w}}{1-|p|^2} \pd{F}{\bar{w}}.
  \]
  Define complex valued 1-forms
  \begin{align*}
     \theta_0 &= dw - 
                    \left (
                       \frac{w \bar{p}}{1-|p|^2}
                       - z 
                    \right ) \, dp
                    - \frac{\bar{w}}{1-|p|^2} \, d \bar{p} \\
     \omega_0 &= dz - F \, dp \\
     \pi_0    &= dp.
  \end{align*}
  These define a section of a unique $G$ structure,
  and moreover the foliation 
  \[
    \theta=\bar{\theta}=\omega=\bar{\omega} = 0
  \]
  if it is a fibration (which is the case locally)
  will have as base a 4-manifold $M$, so that
  these $\theta_0, \omega_0, \pi_0$ are an
  adapted coframing for a pseudocomplex structure
  on $M$. It is remarkable that this structure
  is dual to the previous one via 
  \[
    (W,Z,P)=(w,p,-z)
  \]
  so that all of these examples are in fact
  examples of 4-parameter families of $J$ 
  curves in the previous example. Moreover,
  the generic 4-parameter family of $J$ curves in
  that example will give rise to an
  example of this kind.

    Locally, up to diffeomorphism, these
  examples are described by the choice of 
  function $F$, depending on 2 real functions
  of 5 real variables. The possible symmetry
  pseudogroups are unknown, but there are
  several homogeneous cases (including
  some depending on at least two parameters).
\end{enumerate}

\section{Kobayashi and Brody hyperbolicity}

An $E$ \emph{line} in a manifold $M$
with a pseudocomplex structure $E \subset \Gro{2}{TM}$
is a nonconstant basic 
parameterized $E$ curve $\mathbb{C} \to E$.
A pseudocomplex 4-manifold is called Brody hyperbolic
if it admits no $E$ lines.

It is clear that in the category
of complex manifolds, a holomorphic
fiber bundle with Brody hyperbolic
base is Brody hyperbolic precisely
if its fibers are. Conversely, for
a complex surface which is a 
holomorphic fiber bundle with
hyperbolic fibers, the base is
hyperbolic precisely when the
surface is Brody hyperbolic.
This is not true for pseudocomplex
manifolds, as we see from the
Darboux integrable
almost complex structure \ref{itm:DrbInt}.

We now define the Kobayashi pseudometric. 
A \emph{graph of disks} is a 
connected Riemann surface with two marked points
and finitely many ordinary double
points, whose finitely many 
irreducible components are disks.
We impose the Poincar{\'e} metric on each disk,
so that it is complete with constant negative
curvature $-1$. This induces a metric
on the entire graph of disks.
We define the \emph{length} 
of a graph of disks to be the distance between
the two marked points.

A \emph{Kobayashi chain} (or simply a \emph{chain}) 
is an $E$ curve
\[
\phi : C \to E
\]
where $C$ is a graph of disks.
We say that this
chain joins $p \circ \phi(z) \in M$ to 
$p \circ \phi(w) \in M$
where $z$ and $w$ are the two marked
points.
Define the \emph{Kobayashi metric} 
\[
\operatorname{dis}_{K} (x,y) 
= \inf \operatorname{length} \left \{ \phi_j \right \}
\]
with the infimum taken over all chains joining
$x$ to $y$. The existence of a chain joining
$x$ to $y$ (if $M$ is connected)
is an elementary application of the continuity 
method. Therefore the Kobayashi pseudometric
is defined, symmetric, and satisfies the triangle
inequality. But it may fail to be positive
between distinct points. Note that we could
restrict to using only trees of disks
instead of graphs of disks, i.e.\
ask that the metric space formed
by our graph of disks be strictly convex.
This would provide the same metric.
Or we could be more generous, and allow
any Riemann surface with nodes and 
two marked points whose universal cover
is a graph of disks---a hyperbolic curve
with at most nodal singularities.

We will say that a pseudocomplex manifold
is \emph{Kobayashi hyperbolic} if the Kobayashi
pseudometric gives positive distance
between distinct points.

As usual, we will say that a Riemann surface 
without boundary is hyperbolic
if it admits a metric of constant
negative curvature in its conformal
class.

\begin{proposition} A pseudocomplex manifold
which is Kobayashi hyperbolic is Brody
hyperbolic.
\end{proposition}
\begin{proof} Suppose that our pseudocomplex
manifold fails to be Brody hyperbolic,
so that we have a basic parameterized
$E$ curve $\phi : \mathbb{C} \to E.$
Take a sequence of concentric disks in $\mathbb{C}$
of arbitrarily large radius. Impose the Poincar{\'e}
metric on each, conformally matching the induced
metric. Then any pair of points of $M$ in the
image of $p \circ \phi$ are joined by Kobayashi
chains of arbitrarily small length.
\end{proof}

 Given a parameterized
$E$ disk $\phi : D \to E$
with $D$ bearing the Poincar{\'e}
metric, we define its \emph{velocity}
to be $(p \circ \phi)'(0) \cdot e \in TM$ 
where $0 \in D$
is any fixed point, and $e$ is any fixed
unit vector $e \in T_0 D$. 
If we holomorphically map the disk to itself,
by dilation, we can slow down the
velocity, scaling it by any $t$
with $0 \le t \le 1$. The question is whether
we can speed it up.

\begin{proposition} A compact and proper
pseudocomplex 4-manifold is Brody hyperbolic
precisely when it is Kobayashi hyperbolic.
\end{proposition}
\begin{proof} The argument is identical
to Brody's in \cite{Brody:1978}. We sketch it.
Suppose that \( E \subset \Gro{2}{TM} \)
is our proper pseudocomplex structure on
a compact manifold $M$.
Take any Riemannian metric on $M$.
Given a sequence of parameterized $E$ disks,
\( \phi_j : D \to E, \) 
with arbitrarily large velocity vectors,
Brody \cite{Brody:1978}
presents an argument that shows that we can
replace these $\phi_j$ with new parameterized
$E$ disks $\psi_j$ so that 
\[
\left ( p \circ \psi_j \right )'
\]
is largest at the origin of the disk.
and so that the norm of the velocity of $\psi_j$
becomes arbitrarily large. We can then 
reparameterize to produce,
just as Brody does, a sequence of parameterized
$E$ disks $\xi_j : D_j \to E$, using disks $D_j$ in 
the complex plane of arbitrarily
large radii, so that the velocities of these $\xi_j$
disks are bounded. Since their first derivatives
are largest at the origin, this ensures that
their first derivatives are bounded, so that
some subsequence converges uniformly. The uniformity
ensures by Corollary \ref{cor:UnifCon}
that they converge uniformly with all derivatives
to an $E$ line \( \xi : \mathbb{C} \to E. \)
\end{proof}

Given an infinitesimal symmetry of a 
pseudocomplex structure \( E \subset \Gro{2}{TM} \)
described by a vector field $X$ on $M$,
we can prolong $X$ to produce a vector
field $\hat X$ on $E$, defined by
the requirement that the flow of
$\hat X$ on $E$ is given by differentiating
the flow of $X$, and applying its
derivative to a 2-plane of tangent vectors: 
\[
\exp \left ( t \hat X \right ) e 
=
\left ( \exp \left ( t X \right ) \right )'(m) \cdot e
\]
for \( e \in E_m \subset \Gro{2}{T_m M}. \)
This flow is defined through each point,
for small values of $t$, not only for real
$t$ values by also for complex $t$ values
since $\hat X$ is easily seen to be 
a $J$ holomorphic vector field.
The flow defines a foliation of $E$ by
$J$ holomorphic curves, wherever \( \hat X \ne 0, \)
in other words above points of $M$ where $X \ne 0$.
The lift $\hat X$ is complete precisely
when $X$ is.

\begin{theorem} The symmetry group of
a Brody hyperbolic 
proper pseudocomplex structure
on a compact manifold is
finite.
\end{theorem}
\begin{proof} Take an infinitesimal symmetry $X.$
The flow of $\hat X$ is
complete, by compactness, for $X$ any
infinitesimal symmetry. Not every flow
curve of $\hat X$ is actually a
parameterized $E$ curve
for $M$. However, if we pick a 
flow curve $\Sigma$ of $X$, passing
through a point where $X \ne 0$, it is easy to
see that it will be an 
$E$ curve, and that
its lift $\hat \Sigma$ 
will be  a flow
curve of $\hat X$. The completeness
ensures that its domain will be
the entire complex line. Therefore $X$
must vanish everywhere.

  By ellipticity of the $G$ structure
on $E$, the group of symmetries of 
a proper pseudocomplex structure
on a compact manifold is a finite-dimensional
Lie group (see \cite{Ochiai:1966}).
But without any nonzero infinitesimal
symmetries, it must be zero dimensional,
so a discrete group. 

By the
uniform convergence theorem for symmetries
(Theorem \ref{prop:CmptSym}), we known 
that a uniform limit of symmetries
converges with all derivatives. Suppose that
we have a sequence of symmetries.
We can bound their derivatives by
taking an $E$ disk, and looking
at how its velocity gets stretched
under a symmetry.
The stretching factor must be bounded,
or we will be able to recursively
apply the symmetry to generate 
$E$ disks of arbitrarily large velocity.
This bounds the derivative of the
symmetry, since there is an $E$ disk
with any velocity close enough to $0$,
by continuity. Therefore any sequence
of symmetries have a convergent 
subsequence, so the symmetry group
is compact. Because it is also discrete,
it must be finite.
\end{proof}

\begin{corollary}
Let $M$ be a 4-manifold and 
\( E \subset \Gro{2}{TM} \) a proper
pseudocomplex structure on $M$.
The universal cover of $M$, with
the induced proper pseudocomplex
structure, is compact and Brody
hyperbolic precisely if $M$ is,
which happens precisely when
both are compact and Kobayashi
hyperbolic.
\end{corollary}

\begin{proposition} A compact and proper 
pseudocomplex 
4-manifold, equipped with an arbitrary
Finsler metric $v \mapsto |v|$, is
Kobayashi hyperbolic precisely when
the norms $|v|$ of velocities of 
$E$ disks are bounded.
\end{proposition}
\begin{proof} The proof consists in an elementary
limit argument. See \cite{Brody:1978}.
\end{proof}

We now define the Royden pseudonorm on $TM$.
Given any vector $v \in TM$, we can also ask whether
it can be scaled to become the velocity
vector of a parameterized $E$ disk. 
The Royden pseudonorm $|v|$ of any vector 
$v \in TM$ is defined to be the largest
scale factor $t \ge 0$ so that $tv$ occurs
as the velocity of a disk. 

  A simple application of the continuity
method (for closed disks) shows that $|v| > 0$ for $v \ne 0$,
so that the problem with this norm is
only that it may be infinite: $|v|=\infty$,
as occurs for any tangent vector in $M=\mathbb{C}^2$.
Define lengths of paths $\gamma : [a,b] \subset \mathbb{R} \to M$
to be
\[
\operatorname{length} \gamma
=
\int_a^b | \dot \gamma(t) | \, dt
\]
and define the Royden 
pseudodistance between points to be
the infimum of lengths of paths joining them.
Another application of continuity
shows that this length coincides
with the Kobayashi pseudodistance; 
see \cite{KruglikovOverholt:1997} for
details.

\section{Bangert's theory of lines in a torus}

Let $\Lambda \subset \R{2n}$ be a lattice.
We will say that a symplectic structure
on the torus $T^{2n} = \R{2n}/\Lambda$ is
\emph{standard} if it pulls back to the
standard $dq^i \wedge dp^i$ structure
on $\R{2n}.$ It is unknown
if there are any nonstandard symplectic
structures.

\begin{theorem} A proper pseudocomplex structure
on $T^{4}$ tamed by a standard symplectic
structure is not Brody hyperbolic.
\end{theorem}
\begin{proof} The proof is the same as in 
\cite{Bangert:1998}, using the quasiminimality of
$E$ curves, the pseudoconvexity of small balls
(see Section \vref{sec:PLUSH}), and Gromov compactness.
\end{proof}
This result shows that while Brody hyperbolic
structures are an open set, they are not always
dense. In fact, it is not very clear if there
are any Brody hyperbolic pseudocomplex structures on the
torus.

It remains a mystery (first considered by J{\"u}rgen Moser) 
whether entire foliations
of complex tori survive (perhaps in a KAM sense)
perturbation of the complex structure to
a pseudocomplex structure.

\section{Plurisubharmonic functions 
and pseudoconvexity of hypersurfaces}
\label{sec:PLUSH}

\subsection{Pseudoconvexity}
To define pseudoconvexity of hypersurfaces $H \subset M$,
we first need to lift $H$ up to $E \to M$. The
tangent space of $H$ at any point $m \in H$ contains
a discrete set of 2-planes belonging to $E_m \subset \Gro{2}{T_mM}$.
We construct $\tilde{H} \subset E$ to be the
set of all of these 2-planes. It is easy
to see that $\tilde{H}$ is a submanifold of $E$
of dimension 3, and $p : \tilde{H} \to H$
given by $p : E \to M$ is a local diffeomorphism. 
Moreover, $\tilde{H}$ is totally real for the
almost complex structure of $E$. If $E$ is
a proper pseudocomplex structure on $M$, then
in fact $\tilde{H} \to H$ is a diffeomorphism.

Along $H$, the $G$ structure $B$ on $E$ can
be reduced to a principal $G_1$ bundle, where
$G_1$ is the set of complex matrices
\[
\begin{pmatrix}
a & 0 & 0 \\
b & c & 0 \\
\frac{a \bar{b}iL}{|c|^2} & 0 & \frac{a}{c}
\end{pmatrix}
\]
where $a$ is real, $b,c$ are complex,
and $L$ is a relative invariant, which
we will call the \emph{Levi invariant}.
Moreover, 
\begin{align*}
\theta &= \bar \theta \\
\pi &= i L \bar \omega  \; .
\end{align*}
This is easy to calculate (see \cite{McKay:1999}
for details). If $L \ne 0$ at all points
of the bundle over $H$, then we will say
that $H$ is strictly pseudoconvex.

Following Cartan's work on
real hypersurfaces in complex surfaces
\cite{Cartan:136, Cartan:136bis}, we
find the following structure equations
on the reduced bundle of any pseudoconvex $H$:
\begin{align*}
d 
\begin{pmatrix}
\theta \\
\omega \\
\beta \\
\gamma \\
\phi
\end{pmatrix}
&=
\begin{pmatrix}
- ( \gamma + \bar \gamma ) \wedge \theta + i \, \omega \wedge \bar \omega \\
- \beta \wedge \theta - \gamma \wedge \omega \\
- \phi \wedge \omega + \bar \gamma \wedge \beta 
  + R \theta \wedge \bar \omega \\
- \phi \wedge \theta - 2i \bar \beta \wedge \omega 
  - i \beta \wedge \bar \omega \\
i \beta \wedge \bar \beta + (\gamma + \bar \gamma ) \wedge \phi
+ (S \omega + \bar S \bar \omega ) \wedge \theta
\end{pmatrix}
\\
dR &= R (\gamma + 3 \bar \gamma) - \bar S \omega - A \theta - B \bar \omega \\
dS &= S (3 \gamma + 2 \bar \gamma) - i \bar{R} \beta 
  - E \omega -F \bar \omega - G \theta
\end{align*} 
where $\pi = i \bar \omega$, and $A,B,E,G,R,S$
are complex valued functions, and $F$ is a real valued
function (on the induced bundle over $H$).
These are the same structure equations one
usually encounters in the theory of CR geometries
for real hypersurfaces in complex surfaces. 

\subsection{Pseudoconvex foliations}

If we have a foliation $F$ by pseudoconvex
hypersurfaces, we can construct a 
submanifold $\tilde{F} \subset E$
consisting of the $\tilde{H}$ submanifolds
where $H$ runs over all of the hypersurfaces
of the foliation. This $\tilde{F}$ is then
a foliated 4-manifold. On each leaf, we can
produce the subbundle of our $G$ structure
described above for hypersurfaces. Then
this produces a subbundle of our $G$ bundle
over all of $\tilde{F}$, and on this
subbundle, $\pi = i \bar \omega$. Moreover,
on each leaf, $\theta = \bar \theta$. 

A \emph{defining function} for a foliation
$F$ by hypersurfaces is a function $f : M \to \mathbb{R}$
so that the level sets of $f$ are the leaves
of $F$. Locally, a defining function exists.

\begin{proposition} A defining function
of a foliation by pseudoconvex hypersurfaces
in a manifold $M$ with pseudocomplex structure
$E$ must pull back to any $E$ curve to have
positive Laplacian, at all nonsingular points of the
curve where the tangent plane to the
projection of the curve to $M$ 
is nearly tangent to a level 
set of the defining function.
\end{proposition}
\begin{proof} Nearly tangent means
that the tangent plane of the $E$
curve should be close to $\tilde{F}$,
where $F$ is the foliation. 
Without loss of generality,
suppose that the $E$ curve
actually
passes through a point of $\tilde{F} \subset E$. 
It is easy to calculate the
Laplacian of $f$ in terms of the 
subbundle structure equations:
\[
- i \bar{\partial} \partial f
= i \omega \wedge \bar \omega 
\mod{\theta}
\]
so that if $z$ is a holomorphic coordinate
on an $E$ curve, and $\omega = g \, dz$, then
\[
\Delta f = \frac{1}{2} |g|^2
\]
on the $E$ curve.
\end{proof}

\begin{theorem} Consider a nonconstant 
parameterized $E$
curve $\phi : C \to E$. Suppose that the
interior of $C$ is mapped by $p \circ \phi : C \to M$
into the closure of a region
$U \subset M$ with smooth 
strictly pseudoconvex boundary.
\[
\xymatrix{
C \ar[r]^{\phi} \ar[dr]_{p \circ \phi} & E \ar[d]^p \\
& U \subset M
}
\]
No interior point of $C$ gets mapped to the
boundary of $U$, and if a boundary point
of $C$ gets mapped to a boundary point of
$U$:
\[
z \in \partial C \mapsto p \circ \phi(z) = m \in \partial U,
\]
then $p \circ \phi$ is an immersion near $z$, and
\[
(p \circ \phi)'(z) \cdot T_z C \not\subset T_m \partial U
\]
\end{theorem}
\begin{proof} As in the complex case.
\end{proof}

\subsection{Plurisubharmonic functions}

A \emph{plurisubharmonic function}
is a function $f : M \to \mathbb{R}$
whose restriction to any nonsingular
$E$ curve has positive Laplacian.
The story of plurisubharmonic
functions is radically different from the
almost complex case.

\begin{theorem} Let $E$ be a proper pseudocomplex
structure on a 4-manifold $M$. There is
a plurisubharmonic function near any
point of $M$ precisely when $M$ is almost
complex.
\end{theorem}
\begin{proof} Take a function $f : M \to \mathbb{R}$,
and pull it back to $E$. Now differentiate it:
\[
df = f_1 \theta + f_{\bar 1} \bar \theta + 
f_2 \omega + f_{\bar 2} \bar \omega
\]
and then take its second derivative:
\[
d \bar \partial f
=
f_{2 \bar 2} \omega \wedge \bar \omega
+
\bar{f_2} \bar{S_2} \bar \pi \wedge \omega 
-
f_2 S_2 \pi \wedge \bar \omega \; .
\]
If $f_2 S_2(e) \ne 0,$ then there is a 2-plane
in $T_e E$ on which $\theta=\bar \theta=0$
and on which 
\[
-i d \bar \partial f < 0
\]
and another on which 
\[
-i d \bar \partial f > 0 \; .
\]
By local solvability of elliptic partial
differential equations, we can then
find $E$ curves on which the Laplacian
takes either sign. Therefore we will
need $f_2 S_2 = 0$ at each point to
have plurisubharmonicity of $f$.
But it is easy to calculate that $f_2=0$
forces $f$ to be constant. Therefore
we need $S_2=0$. 
By Proposition \ref{prop:WhenAC} (assuming
$E$ proper), this occurs precisely for
$E$ almost complex. 

  Conversely, if $E$ is almost complex,
then in adapted coordinates the function
\[
f(z,w) = |z|^2 + |w|^2
\]
is plurisubharmonic near the origin.
\end{proof}

\section{Bishop disks}

Given a real surface $\Sigma \subset M$ immersed
in a 4-manifold $M$ with pseudocomplex
structure $E$, we can define its \emph{lift}
to be $\tilde{\Sigma}$ to be the set of
its tangent planes, thought of as a surface
inside $\Gro{2}{T_mM}$. Under the projection
\[
\Gro{2}{TM} \to M
\]
this $\tilde{\Sigma}$ is taken diffeomorphically
to $\Sigma$. We call a point of intersection of
$\tilde{\Sigma}$ with $E$ an \emph{elliptic
point} if it is a positive intersection,
a \emph{hyperbolic point} if it is a negative
intersection, and otherwise we call it
a \emph{parabolic point}. See \cite{Moser/Webster:1983}
for more information about the theory of
elliptic and hyperbolic points in complex
surfaces.  A point of $\Sigma$ is called
\emph{totally real} if the tangent
space to $\Sigma$ at that point does
not belong to $E$. The surface $\Sigma$
is called a \emph{totally real surface}
if all of its points are totally real.
A point which is not totally real will
be called an \emph{$E$ point} of $\Sigma$.

\subsection{Totally real surfaces}
Suppose that $\Sigma$ is totally real.
We define a submanifold $\hat{\Sigma} \subset E$
to be  the set
of all 2-planes $P \subset T_m M$
which belong to $E$ and strike the tangent
planes of $\Sigma$ in a line. By our
microlocal geometry results, this $\hat{\Sigma}$
is a smooth 3-manifold immersed in $E$.
It is easy to see that $\hat{\Sigma}$
is totally real as a submanifold of the
almost complex manifold $E$.

\begin{proposition} Let $X$ be a 
totally real immersed submanifold
of an almost complex manifold $Z$. 
Every point $x \in X$ has a neighborhood
$U$ in $Z$ in which there are no
pseudoholomorphic curves $\phi : C \to Z$
from a compact Riemann surface $C$ with
boundary, so that 
\(
\phi ( \partial C ) \subset X
\)
and 
\(
\phi (C) \subset U .
\)
\end{proposition}
\begin{proof} Take adapted complex
coordinates 
\[
z^1, \dots, z^n
\]
as in \cite{McDuff:1994},
so that the totally real submanifold
becomes Lagrangian (indeed, we can
ask that it become the set of real
points of $\mathbb{C}^n$). Then
integrate the K{\"a}hler  form
\[
\frac{i}{2} \sum_j dz^j \wedge dz^{\bar j}
\]
which must be positive on a ``small''
pseudoholomorphic curve, because
the curve is nearly holomorphic.
But the boundary sits in a Lagrangian
manifold, which (if we use a small
enough neighborhood) contradicts
the quasi-minimality of the pseudoholomorphic
curve.
\end{proof}

\begin{corollary} The same result is true
for totally real surfaces in pseudocomplex
4-manifolds.
\end{corollary}
\begin{proof} Looking upstairs in $E$,
an $E$ curve with boundary in a totally
real surface $\Sigma \subset M$ becomes a pseudoholomorphic
curve in $E$ with boundary in $\hat{\Sigma}$. 
\end{proof}

\subsection{Elliptic points and Bishop disks}

When $\Sigma$ has an $E$ point,
$\hat{\Sigma}$
is not well-defined. We can define another
notion of lifting $\Sigma$ to $E$:
let $\Sigma' \subset E$ simply be the
pullback bundle of $E \to M$ under
$\Sigma \to M$, i.e.\ $p : E \to M$
is the projection, and $\Sigma' = p^{-1} \Sigma$.

Pulling back the $G$ structure to $\Sigma'$,
one finds that the structure equations give,
near a point $e_0 \in \Sigma'$ which lies
above an $E$ point of $\Sigma$:  
\[
\theta = f \omega + g \bar \omega
\]
where $f, g : B|_{\Sigma'} \to \mathbb{C}$
are smooth functions vanishing at $e_0$.
We say that the point $e_0$ is \emph{nondegenerate}
if
\[
df \wedge d \bar f \wedge dg \wedge d \bar g \ne 0
\]
at $e_0$. Write
\[
df \wedge d \bar f \wedge dg \wedge d \bar g =
h \omega \wedge \bar \omega \wedge \pi \wedge \bar \pi \; .
\]

\begin{proposition} An $E$ point $e_0$ is 
elliptic precisely when $h > 0$ and
hyperbolic precisely when $h < 0$.
\end{proposition}
\begin{proof} This requires only calculating
examples of surfaces with arbitrary 2-jet,  
since it depends only on the 2-jet of any surface.
For more details, see \cite{McKay:1999}.
\end{proof}

At an elliptic point, we can further refine
the structure equations (see \cite{McKay:1999}),
obtaining at the point $e_0$:
\begin{align*}
&d 
\begin{pmatrix}
\omega \\
\pi \\
f \\
g \\
C
\end{pmatrix}
\\
&=
\begin{pmatrix}
- \gamma \wedge \omega - S_2 \pi \wedge \bar \omega \\
- \bar \gamma \wedge \pi + \left (
	2 \bar S_2 + \bar C T_3 \right ) \omega \wedge \bar \pi 
      + \left ( \bar C \bar S_2 + T_3 \right ) \bar \pi \wedge \pi
	+  q \left ( \pi - \bar \omega \right ) \wedge \omega
	+ r \pi \wedge \bar \omega \\
- \pi + \bar \omega \\
\omega + C \bar \omega \\
C (\bar \gamma - \gamma )
+ (r-Cq-T_2) \omega
+(s-Cr) \bar \omega
+ 2S_2 \pi
- C (T_3 + C S_2) \bar \pi
\end{pmatrix}
\end{align*}

\begin{theorem}[Bishop--Ye]
Let $M$ be a 4-manifold with pseudocomplex
structure $E$, and $\Sigma \subset M$ an
immersed surface with elliptic point $m \in \Sigma$.
There is a smooth embedding $\phi : X \to M$ of  a half
ball $X \subset \mathbb{R}^3$ so that $\phi(0) = m$,
and $\phi$ restricts to each half sphere to
be an embedded $E$ curve with boundary in $\Sigma$.
\end{theorem}
\begin{proof} The argument is the same as
in \cite{Ye:1998}. Also see \cite{McKay:1999}
for more details. Essentially the idea is
the same as in \cite{Bishop:1965}: the linearization
of the problem of constructing the half spheres
has one dimensional kernel and zero dimensional
cokernel. This allows perturbing everything
away from the flat case: $M=\C{2},$
\[
\Sigma = \left \{ 
w = |z|^2 + \frac{\lambda}{2} \left ( z^2 + \bar{z}^2 \right )
\right \}
\]
and
\[
X = \left \{
w = t, 
|z|^2 + \frac{\lambda}{2} \left ( z^2 + \bar{z}^2 \right ) < t
\, | \, t \in \R{}
\right \} \; .
\]
We can approximate in adapted coordinates
(as in Corollary \ref{cor:AdaptedCoords})
with this example (see \cite{Moser/Webster:1983}). 
What makes the a priori estimates work is that
very small $E$ disks near to the origin
have small Sobolev norm of the (nonlinear) operator
\[
\sigma[w] = w_{\bar z} - Q \left ( z,w, w_z \right )
\]
while they have Sobolev norm of roughly
fixed magnitude for
\( \sigma'[w], \)
so that we can move $w(z)$ just a little, 
and effectively push toward
$\sigma=0$. 
\end{proof}

\begin{theorem}[Uniqueness of Bishop disks]
There are at most two possible embedded maximal half
spheres of Bishop disks as defined in the 
previous theorem, one on each side of a 
real surface with an elliptic point. 
\end{theorem}
\begin{proof}
To obtain the uniqueness of Bishop disks,
we need to understand the intersections
of such disks. Suppose that we have
two disks with boundary in a surface,
so that both boundaries are smooth curves
and lie close
to an elliptic point. Suppose further
that the disks are embedded, and transverse
to the real surface along their boundaries.
Neither disk
can touch the elliptic point itself
since that would require tangency with
the real surface at that point. Suppose that
our disks have an interior intersection.
Then we know (by McDuff's arguments as
adapted above to our situation) that the
intersection survives small perturbation.
However, we will need to handle boundary
intersections as well. 

For example, consider the real surface
$\mathbb{RP}^2 \subset \mathbb{CP}^2$. Ignore
for the moment the fact that it has no
elliptic point, and think about two 
complex disks with boundary in this real
surface. For example, take $\mathbb{CP}^1 \subset \mathbb{CP}^2$
and slice it into two disks by cutting out
the real points. Then we have two intersecting
closed disks, which may cease to intersect
if we perturb either one slightly.
The delicate part is to make the disks
intersect ``on the same side.''

Let us follow essentially the argument
of Ye \cite{Ye:1998}.
Because the Fredholm
theory tells us that the deformation theory
of these disks is unobstructed, we can make
them both slightly larger, extending them
so that they pass across the real surface.
Then at an intersection point, we can use
adapted coordinates, holomorphic up to leading
order terms, to see that the intersection
of the two disks looks to leading order
just like the complex case. But there is
a further innovation we can introduce here:
we can adapt our coordinates to the real
surface near the intersection point. Since
it is totally real near that point, there
is no obstruction at any finite order to
finding adapted coordinates for which the
totally real surface is just the set of
points at which the adapted coordinates
take on real values. This is easy to
check with the Cartan--K\"ahler theorem.
However, we have
to be very careful: there could be an
obstruction beyond all orders. But we can
pick an order at which to match up the
real surface to the set of real points
of our adapted coordinates, higher
than the order of intersection of our
disks at that point, and thereby pretend
that the real surface is exactly the
set of real points of our coordinates,
with no loss of generality. Then we see
that, with a little manipulation of the
coordinates, one disk looks like
\[
w = a z^k + O \left ( |z|^{k+1} \right )
\]
with $a \in \mathbb{R}$
while the other looks like
\[
w = 0 \; .
\]
The integer $k$ must be at least one.
Under small perturbation through $E$ disks
we see that the
intersection point can at worst break into $k$
intersection points, distributed around
a circle nearly evenly. As in complex geometry,
it is clear then that the intersection
persists under perturbation, and that moreover
the intersection points are nearly complex
conjugates of each other in these coordinates.
In particular, intersections cannot disappear
off the boundary of the disks, since they
will actually ``bounce off'' the boundary.

Therefore if two $E$ disks with boundary
in a totally real surface intersect, then
they continue to intersect after both of them 
receive a small perturbation, as long as the
intersection was at finite order. However,
as our example shows, they might fail to intersect
after perturbation if they start off with infinite order intersection.
By Aronszajn's lemma, this can only occur if they
can be extended to agree, so that we can glue them
together as in the example.

If we had two half balls of Bishop disks, then
either we could glue them together (if they lay
on opposite sides of the real surface) or they
would be identical (near the elliptic point---if
we extend them to be maximal then they would agree)
or one disk from one half sphere family has to have finite
order intersection with one from the other half sphere
family. 
\end{proof}

\section{Directions for further investigations}

The notion that arises from our
work so far is that the moduli spaces
of curves in a symplectic manifold
should be thought of as having curves
on them. While we have only obtained this
result for the unobstructed part of 
four-dimensional moduli spaces 
of curves in a four-dimensional manifold,
it is tempting to believe that this
concept generalizes to arbitrary dimensions.

Recall that the theory of Gromov--Witten invariants
centers on the map from moduli
of curves in a symplectic manifold to moduli
of abstract curves, forgetting the ambient
manifold; (a little) more precisely, if $C \subset X$
is a holomorphic (or $E$) curve in a symplectic
manifold $X$, then let us write $\ModSp{C}{X}$
for the moduli space of deformations of the
curve $C$ as a holomorphic curve in $X$, compactified
by adding stable curves. 
We can also map $C$ to a point, say $*$, so
that we obtain a map
\[
\ModSp{C}{X} \to \ModSp{C}{*}
\] 
by forgetting how $C$ sits in $X$.
This $\ModSp{C}{*}$ is the Deligne--Mumford
compactified moduli space. Gromov--Witten invariants
arise from pulling back cohomology classes:
\[
H^* \left ( \ModSp{C}{*} \right )
\to
H^* \left ( \ModSp{C}{X} \right ).
\]
If, as the author believes, these moduli
spaces $\ModSp{C}{X}$ carry their own
$E$ curves, then we may construct another
moduli space: if $C_1$ is a curve in $\ModSp{C}{X}$,
then we have a space
\[
\ModSp{C_1}{ \ModSp{C}{X} }
\]
and maps
\[
\ModSp{C_1}{ \ModSp{C}{X} } \to \ModSp{C_1}{ \ModSp{C}{*} } 
\]
which generate cohomology classes
\[
H^* \left ( 
\ModSp{C_1}{ \ModSp{C}{*} }
\right )
\to
H^* \left ( 
\ModSp{C_1}{ \ModSp{C}{X} }
\right ) \; .
\]
Presumably this generates a theory like
that of Gromov--Witten invariants.

Another point of view on this story comes
from the theory of partial differential
equations. The duality between the $E$
curves on 4-manifolds explained above
is a (fairly elementary) example of
a B{\"a}cklund transformation. But
the surprise is that it is ``stable.''
Consider the example of $X=\CP{2}$
with the standard complex structure,
and $Y$ the set of lines in $X$.
We can make perturbations
of the Cauchy--Riemann equations
on $X$ to any nonlinear elliptic equations
as long as they are tamed by a symplectic
structure. The space
$Y$ can be replaced by the space of
rational $E$ curves in $X$ in the 
homology class of a complex line,
and the B{\"a}cklund transformation
survives this perturbation. The author
is unaware of any other example
of this sort of stability for
B{\"a}cklund transformations, or of
any evidence that any B{\"a}cklund
transformation is unstable in this
sense.

The author suspects that if we have a dual
pair of compact pseudocomplex four-manifolds, with both
pseudocomplex structures being
proper, then both four-manifolds are 
diffeomorphic to $\CP{2}$.

 It would be nice to know which real analytic
pseudocomplex structures admit
a Schwarz reflection about any totally real
hypersurface, roughly analogous to symmetric spaces.

\nocite{*}
\bibliographystyle{amsplain}
\bibliography{curves}
\end{document}